\documentclass[11pt,subeqn,subarray]{article}
\usepackage{amssymb,amsmath}
\usepackage{graphicx}
\usepackage{flafter}

\newtheorem{theorem}{Theorem}[section]
\newtheorem{corollary}{Corollary}[section]
\newtheorem{lemma}{Lemma}[section]
\newtheorem{definition}{Definition}[section]
\newtheorem{remark}{Remark}[section]

\def\sqr#1#2{{\vcenter{\vbox{\hrule height.#2pt
              \hbox{\vrule width.#2pt height#1pt \kern#1pt \vrule width.#2pt}
              \hrule height.#2pt}}}}
\def\signed #1{{\unskip\nobreak\hfil\penalty50
              \hskip2em\hbox{}\nobreak\hfil#1
              \parfillskip=0pt \finalhyphendemerits=0 \par}}
\def\endpf{\signed {$\sqr69$}}

\def\be{\begin{equation}}
\def\ee{\end{equation}}
\def\bea{\begin{eqnarray}}
\def\eea{\end{eqnarray}}
\def\ca{{\cal A}}

\def\cd{{\cal D}}

\def\ch{{\cal H}}
\def\cl{{\cal L}}

\def\s{\sigma}
\def\l{\lambda}
\def\o{\omega}

\def\L{\Lambda}
\def\D{\Delta}
\def\g{\gamma}
\def\t{\tau}
\def\th{\theta}
\def\m{\mu}
\def\n{\nu}
\def\a{\alpha}
\def\b{\beta}
\def\e{\varepsilon}
\def\f{\varphi}
\def\d{\delta}
\def\z{\zeta}

\def\dbC{\mathbb{C}}

\def\dbE{\mathbb{E}}

\def\dbR{\mathbb{R}}

\def\cA{{\cal A}}

\def\cD{{\cal D}}

\def\cH{{\cal H}}

\def\cK{{\cal K}}
\def\cL{{\cal L}}

\def\D{\Delta}

\def\L{\Lambda}
\def\Si{\Sigma}
\def\F{\Phi}

\def\no{\noindent}

\def\ss{\smallskip}
\def\ms{\medskip}
\def\bs{\bigskip}
\def\hb{\hbox}
\def\q{\quad}
\def\qq{\qquad}
\def\3n{\negthinspace \negthinspace \negthinspace }
\def\2n{\negthinspace \negthinspace }
\def\1n{\negthinspace }

\def\limsup{\mathop{\overline{\rm lim}}}

\def\lan{\mathop{\langle}}
\def\ran{\mathop{\rangle}}

\def\cd{\cdot}
\def\cds{\cdots}

\def\cl{\overline}

\def\Re{{\mathop{\rm Re}\,}}
\def\Im{{\mathop{\rm Im}\,}}

\def\({\Big (}
\def\){\Big )}
\def\[{\Big[}
\def\]{\Big]}
\def\bde{\begin{definition}}
\def\ede{\end{definition}}
\def\be{\begin{equation}}
\def\bel{\begin{equation}\label}
\def\ee{\end{equation}}
\def\bex{\begin{example}}
\def\eex{\end{example}}
\def\bt{\begin{theorem}}
\def\et{\end{theorem}}
\def\bc{\begin{corollary}}
\def\ec{\end{corollary}}
\def\bl{\begin{lemma}}
\def\el{\end{lemma}}
\def\bp{\begin{proposition}}
\def\ep{\end{proposition}}
\def\bas{\begin{assumption}}
\def\eas{\end{assumption}}
\def\br{\begin{remark}}
\def\er{\end{remark}}
\def\ba{\begin{array}}
\def\ea{\end{array}}
\def\ed{\end{document}}
\def\ds{\displaystyle}

\def\ns{\noalign{\ss}}

\begin{document}
\title
{\bf  Regularity Analysis for an Abstract System of Coupled Hyperbolic and Parabolic
Equations }
\author{Jianghao Hao \thanks{School of Mathematical Sciences, Shanxi University,
Taiyuan, Shanxi 030006, China. This author was supported
by NNSFC No. 61374089 and Shanxi Scholarship council of China No.2013-013 .}
\and Zhuangyi Liu \thanks{Department of Mathematics and Statistics,
University of Minnesota, Duluth, MN 55812-2496, USA} \and Jiongmin
Yong \thanks{Department of Mathematics, University of Central
Florida, Orlando, FL 32816, USA. This author was partially supported
by NSF Grant 1007514.}} \maketitle
\date{}

{\small  {\bf Abstract} \vskip2mm In this paper, we provide a
complete regularity analysis for the following abstract system of
coupled hyperbolic and parabolic equations
$$\left\{\2n\ba{ll}
\ns\ds u_{tt}=-Au+\g A^\a w,\\
\ns\ds w_t=-\g A^\a u_t-kA^\b w,\\
\ns\ds u(0)=u_0,\quad u_t(0)=v_0, \quad w(0)=w_0, \ea\right.$$
where $A$ is a self-adjoint, positive definite operator on a complex
Hilbert space $H$, and $(\a,\b)\in[0,1]\times[0,1]$. We are able to
decompose the unit square of the parameter $(\a,\b)$ into three
parts where the semigroup associated with the system is analytic, of
specific order {\it Gevrey} classes, and non-smoothing,
respectively. Moreover, we will show that the orders of {\it Gevrey}
class are sharp, under proper conditions. }

\vskip 2mm


\ms


\section{Introduction}
\setcounter{equation}{0}
\setcounter{theorem}{0}

Let $H$ be a complex Hilbert space with the inner product
$\lan\cd\,,\cd\ran$ and the induced norm $\|\cd\|$. We consider the
following abstract system of coupled hyperbolic and parabolic
equations:
\bel{1.1}\left\{\2n\ba{ll}
\ns\ds u_{tt}=-Au+\g A^\a w,\\
\ns\ds w_t=-\g A^\a u_t-kA^\b w,\\
\ns\ds u(0)=u_0,\q u_t(0)=v_0,\q w(0)=w_0,\ea\right.\ee
where $A$ is a self-adjoint, positive definite (unbounded) operator
on a complex Hilbert space $H$; $\g\ne0$, $k>0$, and $\a,\b\in[0,1]$
are fixed real numbers. Our main interest is the regularity of the
solution to this system in terms of the parameters $\a,\b$.

\ms

We define
$$\cH=\cD(A^{1\over2})\times H\times H.$$
Any element in $\cH$ is denoted by $U=(u,v,w)^T$. Introduce
$$\lan U_1,U_2\ran{}\1n_\cH=\lan A^{1\over2}u_1,A^{1\over2}u_2\ran+\lan
v_1,v_2\ran+\lan w_1,w_2\ran,\qq\forall \ U_i=\begin{pmatrix}u_i\\ v_i\\
w_i\end{pmatrix}\in\cH,~i=1,2.$$
Then $\lan\cd\,,\cd\ran{}\1n_\cH$ is an inner product under which
$\cH$ is a Hilbert space. By denoting $v=u_t$ and
$U_0=(u_0,v_0,w_0)^T$, system (\ref{1.1}) can be written as an
abstract linear evolution equation on the space $\cH$,
\bel{1.2}\left\{\2n\ba{ll}
\ns\ds\frac{dU(t)}{dt}=\cA_{\a,\b}U(t),\qq t\ge0, \\
\ns\ds U(0)=U_0,\ea\right.\ee
where the operator $\cA_{\a,\b}:\cD(\cA_{\a,\b})\subseteq\cH\to\cH$
is defined by
\bel{1.3}\cA_{\a,\b}=\begin{pmatrix}
0 & I & 0\\
-A & 0 & \g A^\a \\
0 & -\g A^\a & -kA^\b \end{pmatrix},\ee
with the domain
\bel{1.4}
\cD(\cA_{\a,\b})=\cD(A)\times\cD(A^{\a\vee{1\over2}})\times
\cD(A^{\a\vee\b}),\ee
where $a\vee b=\max\{a,b\}$ for any $a,b\in\dbR$. It is known that
$\ca_{\a,\b}$ (which is identified with its closure) generates a
$C_0$-semigroup $e^{\ca_{\a,\b}t}$ of contractions on $\cH$
(\cite{ABB}). Then the solution to the evolution equation
(\ref{1.2}) admits the following representation:
$$U(t)=e^{\ca_{\a,\b}t}U_0,\qq t\ge0,$$
which leads to the well-posedness of (\ref{1.2}). With this in hand,
regularity and stability are the most interesting properties for the
solutions to evolution equations that attract people's attention.
Before going further, let us recall some definitions relevant to the
regularity and stability of $C_0$-semigroups.

\ms

{\bf Definition 1.1.} Let $e^{\mathcal{A}t}$ be a $C_0$-semigroup on
a Hilbert space $\cH$.

\ms

(i) Semigroup $e^{\cA t}$ is said to be {\it analytic} if there
exists an extension $T(\t)$ of $e^{\cA t}$ to the following set
$$\Sigma_\th\equiv\{\t\in\dbC\bigm||\arg\t|<\th\}\cup\{0\},$$
for some $\th\in(0,{\pi\over2})$ so that for any $x\in\cH$,
$\t\mapsto T(\t)x$ is continuous on $\Sigma_\th$ satisfying the
following semigroup property
$$T(\t_1+\t_2)=T(\t_1)T(\t_2),\qq\forall\t_1,\t_2\in\Sigma_\th,~\hb{with $\t_1+\t_2\in\Sigma_\th$},$$
and $\t\mapsto T(\t)$ is analytic over $\Sigma_\th\setminus\{0\}$ in
the uniform operator topology of $\cL(\cH)$ (the space of all linear
bounded operators from $\cH$ to $\cH$).

\ms

(ii) Semigroup $e^{\cA t}$ is said to be of {\it Gevrey} class $\d$
(with $\d>1$) if it is infinitely differentiable and for any compact
set $\cK\subset(0,\infty)$ and any $\th>0$, there exists a constant
$K=K(\th,\cK)$, such that
\bel{1.5}\|\cA^ne^{\cA t}\|_{\cL(\cH)}\le K\th^n(n!)^\d,\qq\forall
t\in\cK,\ n\ge0.\ee

(iii) Semigroup $e^{\cA t}$ is said to be {\it differentiable} if
for any $x\in\cH$, $t\mapsto e^{\cA t}x$ is differentiable on
$(0,\infty)$.

\ms

(iv) Semigroup $e^{\cA t}$ is said to be {\it exponentially stable}
with decay rate $\o>0$ if there exists a constant $M\ge1$ such that
$$\|e^{\cA t}\|\le Me^{-\o t},\qq t\ge0.$$

\ms

(v) Semigroup $e^{\cA t}$ is said to be {\it polynomially stable of
order $j>0$} if there exists a constant $M>0$ such that
$$\|e^{\cA t}\cA^{-1}\|\le Mt^{-j},\qq t>0.$$

\ms

In the above, the first three notions are about the regularity of
$C_0$-semigroups and the last two notions are about the
asymptotically stability of $C_0$-semigroups. We will see shortly
that these notions are intrinsically related. Note that in
(\ref{1.5}), if $\d=1$, then the semigroup is analytic.

\ms

We now briefly recall some history. In 1981, Chen--Russell
(\cite{CR}) considered the abstract elastic system with direct
damping (the so-called linear oscillator) of following form:
\bel{}{d\over dt}\begin{pmatrix}u\\
v\end{pmatrix}=\cA_\a\begin{pmatrix}u\\ v
\end{pmatrix}=\begin{pmatrix}0&I\\-A&-B_\a\end{pmatrix}
\begin{pmatrix}u\\ v\end{pmatrix}\label{1.6}\ee
on $\ch=D(A^\frac{1}{2})\times H$, where both $A$ and $B_\a$ are
(unbounded) positive definite on a Hilbert space $H$. Two
conjectures for the analyticity of the associated $C_0$-semigroup
$e^{\cA_\a t}$ were posed. It was shown by Huang \cite{H1, H2} and
Huang--Liu \cite{HuL} that if $B_\a$ is equivalent to $A^\a$ (in a
certain sense) with ${1\over2}\le\a\le1$, the semigroup $e^{\cA_\a
t}$ is analytic and exponentially stable. Complete regularity
results for such a system were obtained by Chen--Triggiani
(\cite{CT1,CT2}), which says: When $B_\a$ is equivalent to $A^\a$
(in a certain sense), the associated $C_0$-semigroup $e^{\cA_\a t}$
is analytic for ${1\over2}\le\a\le 1$, is of Gevrey class
$\delta>{1\over2\a}$ for $0<\a<{1\over2}$.

\ms

Having the complete results for system (\ref{1.6}), people naturally
turned the attention to thermoelastic equations, such as string,
beam and plate, and so on. In the early 1990's, Russell \cite{R}
proposed an abstract system of a second order conservative equation
coupled with a first order dissipative equation:
\bel{1.7}\frac{d}{dt}\begin{pmatrix}u\\ v\\ w\end{pmatrix}=
\cA\begin{pmatrix}u \\ v\\ w
\end{pmatrix}=\begin{pmatrix}
0 & I & 0\\
-A & 0 & B \\
0 & -B^* & -D \end{pmatrix}\begin{pmatrix}u \\ v\\ w
\end{pmatrix}.\ee
This can be regarded as a system with indirect damping and velocity
coupling. He pointed out that it is desirable to obtain complete
results for system (\ref{1.7}) similar to the known results for
system (\ref{1.6}). This has motivated studies of system (\ref{1.7})
and/or (\ref{1.1}) since then. For (\ref{1.1}), a complete stability
analysis was carried out by the first two authors of the current
paper in 2013 (see \cite{HL}). To state the result, let us introduce
the following sets which give a partition of the unit square
$[0,1]\times[0,1]$:
\bel{S}\left\{\2n\ba{ll}
\ns\ds S=\Big\{(\a,\b)\in[0,1]\times[0,1]\bigm||2\a-1|\le\b\le 2\a
\Big\},\\
\ns\ds S_1=\Big\{(\a,\b)\in[0,1]\times[0,1]\bigm|
2\a\vee{1\over2}<\b
\Big\},\\
\ns\ds
S_2=\Big\{(\a,\b)\in[0,1]\times[0,1]\bigm|\b<1-2\a,~\b\le{1\over2}\Big\},\\
\ns\ds
S_3=\Big\{(\a,\b)\in[0,1]\times[0,1]\bigm|\b<2\a-1\Big\},\ea\right.\ee
where $a\land b=\min\{a,b\}$, and we recall that $a\vee
b=\max\{a,b\}$. See Figure 1. Note that
$$\big[0,{1\over4}\big)\times\big\{{1\over2}\big\}\subseteq S_2.$$

\begin{figure}[h]
\centering
\includegraphics[scale=1]{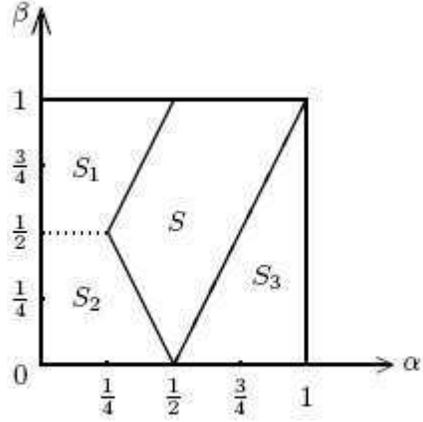}
\caption{Region of stability}
\end{figure}

\ms

\vfil\eject

Here is the stability result found in \cite{HL}.

\ms

\bf Theorem 1.1. \sl The semigroup $e^{\cA_{\a,\b}t}$ has the
following stability properties:

\ms

{\rm(i)} In $S$, it is exponentially stable;

\ms

{\rm(ii)} In $S_1\cup S_2$, it is polynomially stable of order
${1\over2(\b-2\a)}\land{1\over2-2(2\a+\b)}$;

\ms

{\rm(iii)} In $S_3$, it is not asymptotically stable.

\ms

\rm

Note that
$${1\over2(\b-2\a)}\land{1\over2-2(2\a+\b)}=\left\{\2n\ba{ll}
\ns\ds{1\over2(\b-2\a)}>0,\qq(\a,\b)\in
S_1,\\
\ns\ds{1\over2-2(2\a+\b)}>0,\qq(\a,\b)\in S_2.\ea\right.$$

\ms

For the regularity of the semigroup $e^{\cA_{\a,\b}t}$, we recall
the following results from the literature.

\ms

$\bullet$ In 1996, Mu\~{n}oz Rivera and Racke studied the smoothing
property of the semigroup $e^{\cA_{\a,\b}t}$ (\cite{MR}). It was
shown that this semigroup is $C^\infty$ in the region
\bel{S0}S^o=\Big\{(\a,\b)\in
[0,1]\times[0,1]\bigm||1-2\a|<\b<2\a\Big\}.\ee
See Figure 2 in which $S^o$ is shadowed, whose closure is $S$
defined in (\ref{S}).

\begin{figure}[h]
\centering
\includegraphics[scale=1]{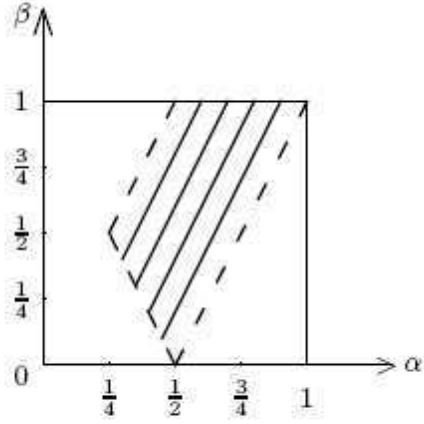}
\caption{Region of $C^\infty$ smoothness}
\end{figure}

\vfil\eject

Now, we divide the unit square $[0,1]\times[0,1]$ further as shown
in Figure 3, where
\bel{R}\left\{\2n\ba{ll}
\ns\ds
R_1=\Big\{(\a,\b)\in[0,1]\times[0,1]\bigm|\a\le\b\le2\alpha-{1\over
2}\Big\},\\
\ns\ds
R_2=\Big\{(\a,\b)\in[0,1]\times[0,1]\bigm|\(2\a-{1\over2}\)\vee{1\over2}<\b<2\a
\Big\},\\
\ns\ds
R_3=\Big\{(\a,\b)\in[0,1]\times[0,1]\bigm|0\le1-2\a<\b\le{1\over2},~
(\a,\b)\ne\({1\over2},{1\over2}\)\Big\},\\
\ns\ds R_4=\Big\{(\a,\b)\in[0,1]\times[0,1]\bigm|0<2\a-1\le\b<\a\Big\},\\
\ns\ds R_5=\Big\{(\a,\b)\in[0,1]\times[0,1]\bigm|0<\b<2\a-1\Big\},\\
\ns\ds R_6=([0,1]\times[0,1])\setminus (R_1\cup R_2\cup R_3\cup
R_4\cup R_5)=S_1\cup S_2 \cup S_I,\ea\right.\ee
with $S_I = ({1\over 2}, 1]\times \{0\}$.
\begin{figure}[h]
\centering
\includegraphics[scale=1]{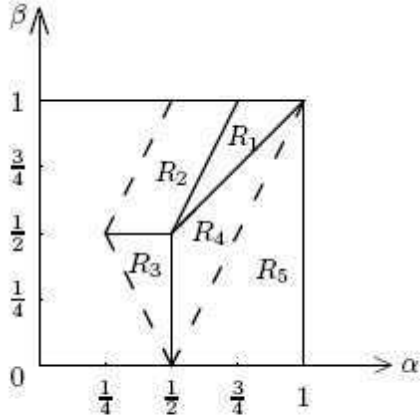}
\caption{Region of regularity}
\end{figure}

\no We see that
$$S^o=R_1\cup R_2\cup R_3\cup R_4,\qq S_3=R_5\cup S_I.$$

$\bullet$ In 1998, Liu and Yong obtained several regularity results
for a general coupled system (\cite{LY}), which implies that the
semigroup $e^{\ca_{\a,\b}t}$ is analytic in $R_1$, and is of {\it
Gevrey} class $\d>\frac{1}{2(2\a-\b)}$ in $R_2$.

\ms

$\bullet$ In 2006, Denk and Racke showed that in region $R_1$ the
semigroup remains analytic in Banach space $L^p(R^n)$, for all
$1<p<\infty$, with $A$ being $-\D$ (\cite{DR}).

\ms

It is natural to ask what can we say about the regularity of the
semigroup $e^{\cA_{\a,\b}t}$ for all the values of $\a,\b\in[0,1]$,
beyond just being analytic in $R_1$ and being $C^\infty$ in $S^o$?
The main results of this paper can be stated as follows.

\ms

\bf Theorem 1.2. \sl The semigroup $e^{\cA_{\a,\b}t}$ has the
following regularity properties:

\ms

{\rm(i)} It is analytic in $R_1$;

\ms

{\rm(ii)} It is of Gevrey class $\d>{1\over\m(\a,\b)}$ in $R_2\cup
R_3\cup R_4\cup R_5$ with
\bel{mu}\m(\a,\b)=\left\{\2n\ba{ll}
\ns\ds2\[(2\a-\b)\land(2\a+\b-1)\],\qq\qq&(\a,\b)\in R_2\cup R_3,\\
\ns\ds{\b\over\a}\,,\qq\qq\qq\q&(\a,\b)\in R_4\cup R_5;\ea\right.\ee

\ms

{\rm(iii)} It is not differentiable in $R_6$.

\ms

\no Moreover, if $A$ admits a sequence of eigenvalues $\m_n\in\dbR$
such that
\bel{1.8}\lim_{n\to\infty}\m_n=\infty,\ee
then the Gevrey class orders in {\rm(ii)} are sharp in the following
sense: For any $\e>0$, the semigroup is not Gevrey class of order
${1\over\m(\a,\b)+\e}$.

\rm

\bs

The significance of the above result includes the following:

\ms

$\bullet$ In the region $R_2\cup R_3\cup R_4$, we establish that
$e^{\cA_{\a,\b}t}$ is of proper order Gevrey classes, instead of
just saying that it is $C^\infty$ as in \cite{MR}.

\ms

$\bullet$ The semigroup $e^{\cA_{\a,\b}t}$ is also shown to be
Gevrey class of a proper order in $R_5$ and not even differentiable
in $R_6$, where, to our best knowledge, there is no
regularity result for the semigroup in the region $R_5\cup R_6$ so far.

\ms

$\bullet$ The Gevrey class orders are sharp for the case that $A$ is
a positive definite self-adjoint operator having a sequence of
(real) eigenvalues that goes to infinite. This is the case when $A$
is a usual elliptic differential operator, say, $-\D$ in a bounded
domain.

\ms

Note that
$${1\over\m(\a,\b)}=\left\{\2n\ba{ll}
\ns\ds{1\over2(2\a-\b)},\qq\q(\a,\b)\in R_2,\\
\ns\ds{1\over2(2\a+\b)-2},\qq(\a,\b)\in R_3.\ea\right.$$
In a word, our results complete the regularity analysis on the
semigroup $e^{\cA_{\a,\b}t}$, in a certain sense. Combining our
results with those found in the literature, we have the following
summary:
\begin{center}
  \begin{tabular}{| c || c | c | }
    \hline
    Regions & Regularity & Stability \\ \hline \hline
    $R_1$ & analytic & exponentially stable \\ \hline
    $R_2$ & Gevrey class $\d>{1\over2(2\a-\b)}$ & exponentially stable \\ \hline
    $R_3$ & Gevrey class $\d>{1\over2(2\a+\b)-2}$& exponentially Stable \\ \hline
    $R_4$ & Gevrey class $\d>{\a\over\b}$ & exponentially stable \\ \hline
    $R_5$ & Gevrey class $\d>{\a\over\b}$& not asymptotically stable \\ \hline
    $S_I$ & not differentiable   &  not asymptotically stable \\ \hline
    $\cl{(S_1\cup S_2)}\cap S$ & not differentiable & exponentially stable \\ \hline
    $S_1$ & not differentiable & polynomially stable of order ${1\over2(\b-2\a)}$ \\ \hline
    $S_2$ & not differentiable & polynomially stable of order ${1\over2-2(2\a+\b)}$ \\ \hline
          \end{tabular}
\end{center}

\ms

The rest of the paper is organized as following. Sections 2 and 3
are devoted to the proof of (i)--(iii) of Theorem 1.2, showing that
the semigroup $e^{\cA_{\a,\b}t}$ has proper regularity in different
regions. Section 4 is for the asymptotic analysis on an eigenvalue
sequence of $\cA_{\a,\b}$, assuming that $A$ has a sequence of
eigenvalues satisfying (\ref{1.8}). Such an analysis will enable us
to show that the orders of Gevrey class obtained in Sections 2 and 3
in different parts of the unit square are sharp.

\ms

\section{Proof of the Main Result}
\setcounter{equation}{0}
\setcounter{theorem}{0}

For the simplicity of presentation, we will take $\g=k=1$ throughout
the rest of the paper.

\ms

In this section, we will present a proof for part (i)--(iii) of
Theorem 1.2. To this end, let us first recall the following standard
result which is stated in a comparable way (see \cite{P,LZ} for
parts (i)--(ii), \cite{T} for part (iii), \cite{LY} for (iv), and
\cite{Borichev-Tomilov 2010} for (v)).

\ms

\bf Lemma 2.1. \sl  Let $\cA:\cD(\cA)\subseteq\cH\to\cH$ generate a
$C_0$-semigroup $e^{\cA t}$ on $\cH$ such that
\bel{2.1}\|e^{\cA t}\|\le M,\qq\forall t\ge0,\ee
for some $M\ge1$ and
\bel{2.2}i\l\in\rho(\cA),\qq\forall\l\in\dbR,~|\l|\hb{ \rm large
enough}.\ee
Then the following hold:

\ms

{\rm(i)} Semigroup $e^{\cA t}$ is analytic if and only if for some
$a\in\dbR$ and $b,C>0$ such that
\bel{2.3}\rho(\cA)\supseteq\Si(a,b)\equiv\Big\{\l\in\dbC\bigm|\Re\l>a-b|\Im\l|\Big\},\ee
and
\bel{2.4}\|(i\l-\cA)^{-1}\|\le{C\over1+|\l|},\qq\l\in\Si(a,b).\ee
This is the case if and only if
\bel{2.5}\limsup_{\l\in\dbR,\,|\l|\to\infty}|\l|\,\|(i\l-\cA)^{-1}\|<\infty.\ee

{\rm(ii)} Semigroup $e^{\cA t}$ is of {\it Gevrey} class $\d>1$ if
and only if for any $b,\t>0$, there are constants $a\in\dbR$ and
$C>0$ depending on $b,\t,\d$ such that
\bel{2.6}\rho(\cA)\supseteq\Si_b(\d)\equiv\Big\{\l\in\dbC\bigm|\Re\l>a-b|\Im\l|^{1\over\d}\Big\},\ee
and
\bel{2.7}\|(i\l-\cA)^{-1}\|\le
C\(e^{-\t\Re\l}+1\),\qq\forall\l\in\Si_b(\d).\ee
This is the case, in particular, if for some $\m\in(\d^{-1},1)$,
\bel{2.8}\limsup_{\l\in\dbR,\,|\l|\to\infty}|\l|^\m\|(i\l-\cA)^{-1}\|<\infty.
\ee

{\rm(iii)} Semigroup $e^{\cA t}$ is differentiable if and only if
for any $b>0$, there are constants $a_b\in\dbR$ and $C_b>0$ such
that
\bel{2.9}\rho(\cA)\supseteq\Si_b\equiv\Big\{\l\in\dbC\bigm|\Re\l>a_b-b\log|\Im\l|\Big\},\ee
and
\bel{2.10}\|(i\l-\cA)^{-1}\|\le
C_b|\Im\l|,\qq\forall\l\in\Si_b,~\Re\l\le0.\ee
This is the case, in particular, if
\bel{2.11}\limsup_{\l\in\dbR,\,|\l|\to\infty}\log|\l|\|(i\l-\cA)^{-1}\|=0.\ee

{\rm(iv)} {\bf(Gearhart--Pruss)} Semigroup $e^{\cA t}$ is
exponentially stable if and only if
\bel{2.12}i\l\in\rho(\cA),\qq\forall \l\in\dbR,\ee
and
\bel{2.13}\limsup_{\l\in\dbR,|\l|\to\infty}\|(i\l-\cA)^{-1}\|<\infty.\ee

\ms

{\rm(v)} {\bf(Borichev--Tomilov)} Semigroup $e^{\cA t}$ is
polynomially stable of order $j>0$ if and only if $(\ref{2.5})$
holds and
\bel{2.14}\limsup_{|\l|\to\infty}|\l|^{-{1\over
j}}\|(i\l-\cA)^{-1}\|<\infty.\ee

\rm

For notational simplicity, hereafter, we write $i\l-\cA$ instead of
$i\l I-\cA$, omitting $I$. In the above result, the regularity and
stability properties of the semigroup $e^{\cA t}$ are deliberately
related to the spectral/resolvent of the generator $\cA$.
Practically, we will use the limit relations (\ref{2.5}),
(\ref{2.8}) and (\ref{2.11}) to establish the regularity property of
the semigroup, and use the spectrum relations (\ref{2.3}),
(\ref{2.6}) and (\ref{2.9}) to show that the relevant indices are
sharp. The following corollary will be useful below.

\ms

\bf Corollary 2.2. \sl {\rm(i)} Suppose $\s(\cA)$ contains a
sequence $\l_n$ such that
\bel{2.15}\lim_{n\to\infty}\Re\l_n=a,\qq\lim_{n\to\infty}|\l_n|=\infty,\ee
for some $a\in\dbR$, then the semigroup $e^{\cA t}$ is not
differentiable.

\ms

{\rm(ii)} Suppose there exists a sequence $\l_n\in\s(\cA)$ such that
\bel{2.16}\limsup_{n\to\infty}{\Re\l_n\over|\Im\l_n|^{1\over\d}}=0.\ee
Then $e^{\cA t}$ is not of Gevrey class $\d$.

\ms

\it Proof. \rm (i) Suppose $e^{\cA t}$ is differentiable. Then for
any $b>0$, there exists an $a_b\in\dbR$ such that
$$\Re\l_n\le a_b-b\log|\Im\l_n|,\qq n\ge1,$$
since $\l_n\in\s(\cA)$. Letting $n\to\infty$ will lead to a
contradiction. Hence the semigroup $e^{\cA t}$ is not
differentiable.

\ms

(ii) We use part (ii) of Lemma 2.1. Suppose $e^{\cA t}$ is of Gevrey
class $\d>0$, then from (\ref{2.6}), for any $b>0$, there exists an
$a\in\dbR$ such that
$$\Re\l_n\le a-b|\Im\l_n|^{1\over\d},\qq\forall n\ge1,$$
since $\l_n\in\s(\cA)$. Thus,
$$0=\limsup_{n\to\infty}{\Re\l_n\over|\Im\l_n|^{1\over\d}}\le-b,$$
a contradiction. \endpf

\ms

We now state two results whose proof will be carried out in the
following section.

\ms

\bf Theorem 2.3. \sl Let
\bel{mu}\m(\a,\b)=\left\{\2n\ba{ll}
\ns\ds1,\qq\qq\qq\qq&(\a,\b)\in
R_1,\\
\ns\ds2\[(2\a-\b)\land(2\a+\b-1)\],\qq\qq&(\a,\b)\in R_2\cup R_3\cup S_1\cup S_2,\\
\ns\ds{\b\over\a}\,,\qq\qq\qq\q&(\a,\b)\in R_4\cup R_5\cup
S_I.\ea\right.\ee
Then
\bel{2.18a}\limsup_{\l\in\dbR,\,|\l|\to\infty}|\l|^{\m(\a,\b)}\|(i\l-\cA_{\a,\b})^{-1}\|<\infty.\ee

\ms

\rm

\bf Theorem 2.4. \sl Let $A$ admit a sequence of eigenvalues
$\m_n\in\dbR$ such that
$$\lim_{n\to\infty}\m_n=\infty.$$
Then there exists a sequence $\l_n\in\s(\cA_{\a,\b})$ of eigenvalues
of $\cA_{\a,\b}$ such that for any $\e>0$,
\bel{2.19}\lim_{n\to\infty}{\Re\l_n\over|\Im\l_n|^{\m(\a,\b)+\e}}=0,\qq
\forall(\a,\b)\in R_2\cup R_3\cup R_4\cup R_5,\ee
and
\bel{2.20}\lim_{n\to\infty}\Re\l_n=a,\qq\lim_{n\to\infty}|\l_n|=\infty,\qq
\forall(\a,\b)\in R_6.\ee

\ms

\rm

To close this section we present a proof of Theorem 1.2.

\ms

\it Proof of Theorem 1.2. \rm Combining Theorem 2.3 and Lemma 2.1,
we obtain that the semigroup $e^{\cA_{\a,\b}t}$ is analytic in
$R_1$, is Gevrey class of order $\d>{1\over\m(\a,\b)}$ in $R_2\cup
R_3\cup R_4\cup R_5$. Also, in $R_6\equiv S_1\cup S_2\cup S_I$,
(\ref{2.20}) holds. Hence, by Corollary 2.2, the semigroup
$e^{\cA_{\a,\b}t}$ is not differentiable there.

\ms

Next, by (\ref{2.19}) and Corollary 2.2, we see that the Gevrey
class order $\d>{1\over\m(\a,\b)}$ of the semigroup for $(\a,\b)\in
R_2\cup R_3\cup R_4\cup R_5$ is sharp. \endpf

\bs

\rm

\ms

We note that
$$\m(\a,\b)=2(2\a-\b)<0,\qq(\a,\b)\in S_1,$$
and
$$\m(\a,\b)=2(\b+2\a)-2<0,\qq(\a,\b)\in S_2.$$
Thus, the corresponding (\ref{2.18a}) implies that the semigroup
$e^{\cA_{\a,\b}t}$ is polynomially stable with order
${1\over2(\b-2\a)}$ and ${1\over2-2(\b+2\a)}$, respectively. The
above two cases are exactly those found in \cite{HL}.

\ms

\section{Analysis on the Resolvent}
\setcounter{equation}{0}
\setcounter{theorem}{0}

In this section, we will prove Theorem 2.3. It is
technical and lengthy. Let us now make some preparations. First of
all, in our proof, the following interpolation theorem will play a
crucial role.

\ms

\bf Lemma 3.1. \sl Let $A:\cD(A)\subseteq H$ be self-adjoint and
positive definite. Then
\bel{inter1}\|A^px\|\le\|A^qx\|^{p-r\over
q-r}\|A^rx\|^{q-p\over q-r},\qq\forall0\le r\le p\le
q,~x\in\cD(A^q).\ee
In particular, for any $\th\in[0,{1\over2}]$, one has (with $r=0$,
$p=\th$, and $q={1\over2}$)
\bel{inter2}\|A^\th
x\|\le\|A^{1\over2}x\|^{2\th}\|x\|^{1-2\th},\qq\forall
x\in\cD(A^{1\over2}),\ee
and for any $\th\in[{1\over2},1]$ (with $r={1\over2}$, $p=\th$, and
$q=1$)
\bel{inter3}\|A^\th
x\|\le\|Ax\|^{2\th-1}\|A^{1\over2}x\|^{2(1-\th)},\qq\forall
x\in\cD(A).\ee

\rm

\ms

The above result is standard. For reader's convenience, we give a
proof here which is very simple and it just costs us a few lines.

\ms

\it Proof. \rm Since $A$ is self-adjoint and positive definite, it
admits a spectrum decomposition. More precisely, there exists a
family of orthogonal projection operators
$\{\dbE_\l,\,\l\in\s(A)\}$, with $\l\mapsto\dbE_\l$ being
nondecreasing such that
\bel{A}Ax=\int_{\s(A)}\l d\dbE_\l x,\qq\forall x\in\cD(A),\ee
where $\s(A)\subseteq(0,\infty)$ is the spectrum of $A$. Clearly,
for any $\th\in\dbR$,
\bel{A(th)x}A^\th x=\int_{\s(A)}\l^\th d\dbE_\l x,\qq
x\in\cD(A^\th).\ee
Now, for any $0\le r\le p\le q$, $x\in\cD(A^q)$, by H\"older's
inequality, one has
$$\ba{ll}
\ns\ds\|A^px\|^2=\int_{\s(A)}\l^{2p}d\|\dbE_\l x\|^2
\le\(\int_{\s(A)}\l^{2q}d\|\dbE_\l x\|^2\)^{p-r\over
q-r}\(\int_{\s(A)}\l^{2r}d\|\dbE_\l x\|^2\)^{q-p\over
q-r}\\
\ns\ds\qq\q~=\|A^qx\|^{2(p-r)\over q-r}\|A^rx\|^{2(q-p)\over
q-r}.\ea$$
This proves (\ref{inter1}). The two special cases (\ref{inter2}) and
(\ref{inter3}) are clear from (\ref{inter1}). \endpf

\ms

Next, for any $\l\in\dbR$, and any
$U\equiv(u,v,w)^T\in\cD(\cA_{\a,\b})$,
\bel{2.7}(i\l-\cA_{\a,\b})U=\begin{pmatrix}i\l&-I&0\\
A&i\l&-A^\a\\0&A^\a&i\l+A^\b\end{pmatrix}\begin{pmatrix}u\\ v\\
w\end{pmatrix}=\begin{pmatrix}i\l u-v\\ Au+i\l v-A^\a w\\
A^\a v+(i\l+A^\b)w\end{pmatrix}.\ee
Our proof for Theorem 2.3 will be based on a contradiction argument.
Suppose for some given $(\a,\b,\m)\in[0,1]\times[0,1]\times[0,1]$,
without having any specific relations among them, the following is
not true:
$$\limsup_{\l\in\dbR,\,|\l|\to\infty}|\l|^\m\|(i\l-\cA_{\a,\b})^{-1}\|
<\infty.$$
Then there exists a sequence
$\{(\l_n,U_n)\bigm|n\ge1\}\subseteq\dbR\times\cD(\cA_{\a,\b})$ with
$U_n\equiv(u_n,v_n,w_n)^T$, and
\bel{}\left\{\2n\ba{ll}
\ns\ds\lim_{n\to\infty}|\l_n|=\infty,\\
\ns\ds\|U_n\|_\cH^2=\|A^{1\over2}u_n\|^2+\|v_n\|^2+\|w_n\|^2=1,\q
n\ge1,\ea\right.\ee
such that
\bel{2.15}\lim_{n\to\infty}|\l_n|^{-\m}\|(i\l_n-\cA_{\a,\b})
U_n\|_\cH=0,\ee
i.e. (note (\ref{2.7}))
\begin{subequations}
\bel{2.16a}\3n\3n\3n\3n\3n\3n\3n\3n\3n\3n
i\l_n|\l_n|^{-\m}A^{1\over2}u_n-|\l_n|^{-\m}A^{1\over 2}v_n=o(1),\ee
\bel{2.16b}\3n
i\l_n|\l_n|^{-\m}v_n+|\l_n|^{-\m}Au_n-|\l_n|^{-\m}A^{\a}w_n=o(1),\ee
\bel{2.16c} i\l_n|\l_n|^{-\m}w_n+|\l_n|^{-\m}A^\a
v_n+|\l_n|^{-\m}A^\b w_n=o(1).\ee
\end{subequations}
Hereafter $o(1)$ stands for a vector in $H$ (or a quantity in
$\dbR$) which goes to zero as $n\to\infty$. The advantage of using
such a notation is that (\ref{2.16a})--(\ref{2.16c}) can be regarded
as a system of equations, which will be convenient below. For the
sequence $\{(\l_n,u_n,v_n,w_n)\}$ satisfying
(\ref{2.16a})--(\ref{2.16c}), we have the following result.

\ms

\bf Lemma 3.2. \sl The following is true:
\begin{subequations}
\bel{s(1)}\3n\3n\3n\3n\3n\3n\3n\3n\3n\3n\3n\3n\3n\3n\3n\3n\3n\3n\3n\3n\3n\3n\3n\3n\3n\3n\3n\3n
i\l_n|\l_n|^{-\m}\|A^{1\over2}u_n\|^2-|\l_n|^{-\m}\lan
v_n,Au_n\ran=o(1),\ee
\bel{s(2)}\3n\3n\3n\3n\3n\3n\3n\3n\3n\3n\3n\3n\3n
i\l_n|\l_n|^{-\m}\|v_n\|^2+|\l_n|^{-\m}\lan
Au_n,v_n\ran-|\l_n|^{-\m}\lan A^\a w_n,v_n\ran=o(1),\ee
\bel{s(3)}\3n\3n\3n\3n\3n\3n\3n\3n\3n\3n\3n\3n\3n\3n\3n\3n\3n\3n\3n\3n\3n\3n\3n\3n\3n\3n\3n\3n\3n
i\l_n|\l_n|^{-\m}\|w_n\|^2 +|\l_n|^{-\m}\lan A^\a
v_n,w_n\ran=o(1),\ee
\bel{s(4)}\1n\3n\3n\3n\3n\3n\3n\3n\3n\3n\3n\3n\3n\3n\3n\3n\3n\3n\3n\3n\3n\3n\3n\3n\3n\3n\3n\3n\3n\3n\3n\3n\3n\3n\3n\3n\3n\3n\3n\3n\3n\3n\3n\3n\3n\3n\3n\3n
|\l_n|^{-\m}\|A^{\b\over2}w_n\|^2=o(1),\ee
\bel{s(5)}\3n\3n\3n\3n\3n\3n\3n\3n\3n\3n\3n\3n\3n\3n\3n\3n\3n\3n\3n\3n\3n\3n\3n\3n\3n\3n\3n\3n\3n\3n\3n\3n\3n\3n\3n\3n\3n\3n\3n\3n\3n\3n
\|A^{1\over2}u_n\|^2+\|w_n\|^2={1\over2}+o(1),\ee
\bel{s(6)}\3n\3n\3n\3n\3n\3n\3n\3n\3n\3n\3n\3n\3n\3n\3n\3n\3n\3n\3n\3n\3n\3n\3n\3n\3n\3n\3n\3n\3n\3n\3n\3n\3n\3n\3n\3n\3n\3n\3n\3n\3n\3n\3n\3n\3n\3n\3n\3n\3n\3n\3n\3n\3n
\|v_n\|^2={1\over2}+o(1),\ee
\bel{s(7)}|\l_n|^{-1}\|A^{1\over2}v_n\|+|\l_n|^{-1}\|Au_n-A^\a
w_n\big\|+|\l_n|^{-1}\|A^\a v_n+A^\b w_n\|=O(1).\ee
\end{subequations}
Hereafter, $O(1)$ stands for a bounded quantity (uniformly in
$n\ge1$) in $\dbR$.

\ms

\it Proof. \rm By taking inner products of (\ref{2.16a}) with
$A^{1\over2}u_n$ and (\ref{2.16b}) with $v_n$, respectively, we
obtain (\ref{s(1)}) and (\ref{s(2)}). Next, by taking inner product
of (\ref{2.16c}) with $w_n$, we have
\bel{2.18}i\l_n|\l_n|^{-\m}\|w_n\|^2+|\l_n|^{-\m}\|A^{\b\over2}w_n\|^2+|\l_n|^{-\m}\lan
A^\a v_n,w_n\ran=o(1).\ee
Adding the obtained (\ref{s(1)}) and (\ref{s(2)}) to (\ref{2.18}),
one has
$$\ba{ll}
\ns\ds|\l_n|^{-\m}\|A^{\b\over2}w_n\|^2+i\[\l_n|\l_n|^{-\m}
\(\|A^{1\over2}u_n\|^2+\|v_n\|^2+\|w_n\|^2\)\\
\ns\ds\qq\qq\qq\qq\qq+2|\l_n|^{-\m}\(\Im\lan Au_n,v_n\ran
+\Im\lan A^\a v_n,w_n\ran\)\]\\
\ns\ds=|\l_n|^{-\m}\|A^{\b\over2}w_n\|^2+i\[\l_n|\l_n|^{-\m}
+2|\l_n|^{-\m}\(\Im\lan Au_n,v_n\ran+\Im\lan A^\a
v_n,w_n\ran\)\]=o(1).\ea$$
Thus, (\ref{s(4)}) follows. Thanks to this equation, (\ref{s(3)})
follows from (\ref{2.18}).

\ms

On the other hand, by taking conjugate of (\ref{s(1)}) and then
multiplying it by $(-1)$, we have
\bel{2.6c}i\l_n|\l_n|^{-\m}\|A^{1\over2}u_n\|^2+|\l_n|^{-\m}\lan
Au_n,v_n\ran=o(1).\ee
By taking conjugate of (\ref{s(3)}) and then multiplying it by
$(-1)$, we have
\bel{2.10c}i\l_n|\l_n|^{-\m}\|w_n\|^2-|\l_n|^{-\m}\lan A^\a
w_n,v_n\ran=o(1).\ee
Combining (\ref{s(2)}) with (\ref{2.6c})--(\ref{2.10c}), one obtains
\bel{}i\l_n|\l_n|^{-\m}\(\|A^{1\over2}u_n\|^2-\|v_n\|^2+\|w_n\|^2\)=o(1),\ee
leading to
\bel{}\|A^{1\over2}u_n\|^2-\|v_n\|^2+\|w_n\|^2=o(1).\ee
Taking into account $\|U_n\|^2_\cH=1$, we obtain
(\ref{s(5)})--(\ref{s(6)}). Finally, by dividing
(\ref{2.16a})--(\ref{2.16c}) by $\l_n|\l_n|^{-\m}$ (note $\m\le1$),
one has
\begin{subequations}
\bel{2.19a}\3n\3n\3n\3n\3n\3n\3n\2n
iA^{1\over2}u_n-\l_n^{-1}A^{1\over 2}v_n=o(1),\ee
\bel{2.19b}\3n iv_n+\l_n^{-1}Au_n-\l_n^{-1}A^{\a}w_n=o(1),\ee
\bel{2.19c}iw_n+\l_n^{-1}A^\a v_n+\l_n^{-1}A^\b w_n=o(1),\ee
\end{subequations}
which implies (\ref{s(7)}). \endpf

\ms

In what follows, for specific situations, we will end up with
$$\hb{either }\q\|A^{1\over2}u_n\|^2+\|w_n\|^2=o(1),\q\hb{ or }\q\|v_n\|^2=o(1),$$
to lead to a contradiction. Now, we present a detailed proof for
Theorem 2.3.

\ms

\it Proof of Theorem 2.3. \rm The proof for $(\a,\b)\in S_1\cup S_2$
can be found in \cite{HL}. We carry out the proof for the rest parts
of the regions in $[0,1]\times[0,1]$. We divide the proof into
several cases.

\ms

\it Case 1. Let $(\a,\b)\in R_1$, i.e., \rm
\bel{}\a\le\b\le2\a-{1\over2},\qq\m=1.\ee
In this case, (\ref{2.16a})--(\ref{2.16c}) are equivalent to
(\ref{2.19a})--(\ref{2.19c}). Since $\a\le\b$, $A^{\a-\b}$ is
bounded. Applying this bounded operator to (\ref{2.19c}), we get
\bel{}iA^{\a-\b}w_n+\l_n^{-1}A^{2\a-\b}v_n+\l_n^{-1}A^\a
w_n=o(1).\ee
Adding the above to (\ref{2.19b}), we obtain
\bel{2.25*}iv_n+\l_n^{-1}Au_n+iA^{\a-\b}w_n+\l_n^{-1}A^{2\a-\b}v_n=o(1).\ee
Furthermore, taking inner product of the above with $v_n$ yields,
\bel{2.26*}i\|v_n\|^2+\l_n^{-1}\lan
A^{1\over2}u_n,A^{1\over2}v_n\ran+i\lan
A^{\a-\b}w_n,v_n\ran+\l_n^{-1}\|A^{\a-{\b\over 2}}v_n\|^2=o(1).\ee
The first and the third terms in the above are clearly bounded.
Making use of (\ref{s(7)}), we see that the second term in the above
is also bounded. So is the fourth term:
\bel{2.27*}|\l_n|^{-1}\|A^{\a-{\b\over2}}v_n\|^2=O(1).\ee
Thus, noting $\m=1$, and using (\ref{s(4)}), we have
$$\ba{ll}
\ns\ds|\l_n|^{-1}|\lan A^\a v_n,w_n\ran|=|\l_n|^{-1}|\lan
A^{\a-{\b\over 2}}v_n,A^{\b\over 2}w_n
\ran|\\
\ns\ds\qq\qq\qq\qq\q\le\big(|\l_n|^{-{1\over2}}\|A^{\a-{\b\over2}}v_n\|\big)
\big(|\l_n|^{-{\m\over2}}\|A^{\b\over2}w_n\|\big)=o(1).\ea$$
Then (\ref{s(3)}) implies
\bel{wn=0}\|w_n\|=o(1),\ee
and (\ref{s(2)}) becomes
\bel{2.29*}i\|v_n\|^2+\l_n^{-1}\lan Au_n,v_n\ran=o(1).\ee
Also, since $\a\le\b$, (\ref{2.26*}) implies
\bel{2.30*}i\|v_n\|^2+\l_n^{-1}\lan
Au_n,v_n\ran+\l_n^{-1}\|A^{\a-{\b\over2}}v_n\|^2=o(1).\ee
Combining (\ref{2.29*})--(\ref{2.30*}), one gets
\bel{2.32*}\l_n^{-1}\|A^{\a-{\b\over2}}v_n\|^2=o(1),\ee
which improves (\ref{2.27*}). Moreover, since $\a\le\b$, by
(\ref{wn=0}), we may write (\ref{2.25*}) as
\bel{2.25**}iv_n+\l_n^{-1}Au_n+\l_n^{-1}A^{2\a-\b}v_n=o(1).\ee
Further, since ${1\over 2}\le 2\a-\b\le 1$,
$\|A^{1-(2\a-\b)}u_n\|$ is bounded. Taking inner product
(\ref{2.25**}) with $A^{1-(2\a-\b)}u_n$ in $H$, we obtain
$$i\lan v_n,A^{1-(2\a-\b)}u_n\ran
+\lan\l_n^{-1}A^{1\over2}v_n,A^{1\over2}u_n\ran+|\l_n|^{-1}\|A^{1-(\a-{\b\over
2})}u_n\|^2=o(1).$$
The first two terms in the above are bounded. So is the third term.
Therefore, making use of (\ref{2.32*}), we finally obtain
$$|\l_n|^{-1}|\lan Au_n,v_n\ran|=\big(|\l_n|^{-{1\over2}}\|A^{1-(\a-{\b\over2})}u_n\|\big)
\big(|\l_n|^{-{1\over2}}\|A^{\a-{\b\over2}}v_n\|\big)=o(1).$$
Then (\ref{2.29*}) implies
\bel{}\|v_n\|^2=o(1),\ee
which is a contradiction to (\ref{s(6)}).

\ms

\it Case 2. Let $(\a,\b)\in R_2$, i.e., \rm
\bel{2.35*}\(2\a-{1\over2}\)\vee{1\over2}<\b<2\a,\qq\m\equiv\m(\a,\b)=2(2\a-\b).\ee
Note that (see Fig.2) in the current case,
\bel{a<b}\a<\b.\ee
From (\ref{2.35*}), one has
$$0<{\m\over4}=\a-{\b\over2}={1\over2}-{1\over2}(\b-2\a+1)<{1\over2}.$$
Thus, by interpolation, using (\ref{s(4)}) and (\ref{s(7)}), we have
$$\ba{ll}
\ns\ds|\l_n|^{-\m}|\lan A^\a
v_n,w_n\ran|\le|\l_n|^{-{\m\over2}}\|A^{\a-{\b\over2}}v_n\|\big(|\l_n|^{-{\m\over2}}
\|A^{\b\over2}w_n\|\big)\\
\ns\ds\le|\l_n|^{-(2\a-\b)}\|A^{1\over2}v_n\|^{2\a-\b}\|v_n\|^{1-2\a+\b}o(1)\le\big(|\l_n|^{-1}\|A^{1\over2}v_n\|\big)^{2\a-\b}o(1)=o(1).\ea$$
Consequently, (\ref{s(3)}) can be written as
\bel{w=0*}\l_n|\l_n|^{-\m}\|w_n\|^2=o(1),\ee
and (\ref{s(2)}) can be written as
\bel{v=0*}i\l_n|\l_n|^{-\m}\|v_n\|^2+|\l_n|^{-\m}\lan
Au_n,v_n\ran=o(1),\ee
which implies
\bel{v=0**}i\|v_n\|^2+\l_n^{-1}\lan Au_n,v_n\ran=o(1).
\ee
We now show that
\bel{Auv=0}\l_n^{-1}\lan A u_n,v_n\ran=o(1).\ee
Since $\a<\b$, $A^{\a-\b}$ is bounded. Applying $A^{\a-\b}$ to
(\ref{2.16c}), we have
$$i\l_n|\l_n|^{-\m}A^{\a-\b}w_n+|\l_n|^{-\m}A^{2\a-\b}v_n+|\l_n|^{-\m}A^\a w_n=o(1).$$
Adding the above to (\ref{2.16b}), one has
\bel{2.33*}i\l_n|\l_n|^{-\m}v_n+|\l_n|^{-\m}Au_n
+|\l_n|^{-\m}A^{2\a-\b}v_n+i\l_n|\l_n|^{-\m}A^{\a-\b}w_n=o(1),\ee
which implies
\bel{2.33**}iv_n+\l_n^{-1}Au_n+\l_n^{-1}A^{2\a-\b}v_n+iA^{\a-\b}w_n=o(1).\ee
Further, by (\ref{w=0*}), $\|w_n\|=o(1)$. Thus, the above becomes
\bel{2.33***}iv_n+\l_n^{-1}Au_n+\l_n^{-1}A^{2\a-\b}v_n=o(1).\ee
By (\ref{s(7)}), for any $0\le\n\le{1\over2}$,
\bel{2.40*}|\l_n|^{-2\n}\|A^\n
v_n\|\le|\l_n|^{-2\n}\|A^{1\over2}v_n\|^{2\n}\|v_n\|^{1-2\n}\le\big(|\l_n|^{-1}\|A^{1\over2}v_n\|\big)^{2\n}
=O(1).\ee
Since $2\a-\b={\m\over2}<{1\over2}$, one has
$$|\l_n|^{-1}\|A^{2\a-\b}v_n\|=|\l_n|^{\m-1}\big(|\l_n|^{-\m}\|A^{\m\over2}v_n\|\big)
=|\l_n|^{-(1-\m)}O(1)=o(1).$$
Thus, (\ref{2.33***}) becomes
\bel{2.33****}iv_n+\l_n^{-1}Au_n=o(1).\ee
Consequently, we obtain
\bel{2.41*}|\l_n|^{-1}\|Au_n\|=O(1).\ee
By interpolation, together with (\ref{2.41*})
$$\ba{ll}
\ns\ds|\l_n|^{-{1+\m\over 2}}\|A^{3+\m\over
4}u_n\|=|\l_n|^{-{1+\m\over 2}}\|A^{1+\m\over
4}(A^{1\over2}u_n)\|\le\|\l_n^{-1}Au_n\|^{1+\m\over 2}\|A^{1\over
2}u_n\|^{1-\m\over 2} =O(1).\ea$$
Now, taking inner product of (\ref{2.33*}) with $|\l_n|^{-{1+\m\over
2}}A^{1+\m\over 4}v_n$ in $H$ leads to
\bel{2.42*}\ba{ll}
\ns\ds i\l_n|\l_n|^{-{1+3\m\over2}}\|A^{1+\m\over
8}v_n\|^2+|\l_n|^{-{1+3\mu\over 2}}\lan
Au_n,A^{1+\m\over4}v_n\ran\\
\ns\ds+|\l_n|^{-{1+3\m\over2}}\|A^{1+3\m\over
8}v_n\|^2+i\l_n|\l_n|^{-{1+3\m\over2}}\lan
A^{\a-\b}w_n,A^{1+\m\over4}v_n\ran=o(1).\ea\ee
Recall that $\b\in({1\over2},1)$ and $\m=2(2\a-\b)$. Thus,
$$\a-\b+{1+\m\over4}={\m\over4}-{\b\over2}+{1+\m\over4}={2\m+1-2\b\over4}
<{\m\over2}<{1\over2},$$
and
$${2-\m-2\b\over2}={1-\m\over2}+{1-2\b\over2}<{1-\m\over2}.$$
Hence, it follow from (\ref{w=0*}) that
$$\ba{ll}
\ns\ds|\l_n|^{{1-3\m\over2}}|\lan
A^{\a-\b}w_n,A^{1+\m\over4}v_n\ran|\2n\2n&=|\l_n|^{2-\m-2\b\over2}\|w_n\|\big(|\l_n|^{-{2\m+1-2
\b\over2}} \|A^{{2\m+1-2\b\over4}}v_n\|\big)\\
\ns\ds&\le|\l_n|^{1-\m\over2}\|w_n\|O(1)=o(1).\ea$$
Then (\ref{2.42*}) becomes
\bel{2.42**}\ba{ll}
\ns\ds i\l_n|\l_n|^{-{1+3\m\over2}}\|A^{1+\m\over
8}v_n\|^2+|\l_n|^{-{1+3\mu\over 2}}\lan
Au_n,A^{1+\m\over4}v_n\ran+|\l_n|^{-{1+3\m\over2}}\|A^{1+3\m\over
8}v_n\|^2=o(1).\ea\ee
Now, taking inner product of (\ref{2.16a}) with $|\l_n|^{-{1+\m\over
2}}A^{3+\mu\over 4}u_n$, we have
\bel{2.43*}\ba{ll}
\ns\ds o(1)=i\l_n|\l_n|^{-{1+3\m\over2}}\lan
A^{1\over2}u_n,A^{3+4\m\over4}u_n\ran-|\l_n|^{-{1+3\m\over2}}\lan
A^{1\over2}v_n,
A^{3+\m\over4}u_n\ran\\
\ns\ds\qq=i\l_n|\l_n|^{-{1+3\m\over2}}\,\|A^{5+\m\over
8}u_n\|^2-|\l_n|^{1+3\m\over2}\lan A^{1+\m\over 4}v_n,
Au_n\ran.\ea\ee
Adding (\ref{2.42**}) to (\ref{2.43*}) and taking its real part, we
get
\bel{}|\l_n|^{-{1+3\m\over2}}\|A^{1+3\m\over 8}v_n\|^2=o(1).\ee
Consequently,
$$\ba{ll}
\ns\ds|\l_n|^{-1}|\lan
Au_n,v_n\ran|\le\big(|\l_n|^{-{3(1-\m)\over4}}\|A^{7-3\m\over8}u_n\|\big)
\big(|\l_n|^{-{1+3\m\over4}}\|A^{1+3\m\over8}v_n\|\big)\\
\ns\ds=|\l_n|^{-{3(1-\m)\over4}}\|A^{3(1-\m)\over8}(A^{1\over2}u_n)\|o(1)
\le\big(|\l_n|^{-1}\|Au_n\|\big)^{{3(1-\m)\over4}}
\|A^{1\over2}u_n\|^{1-{3(1-\m)\over4}}o(1)=o(1).\ea$$
Thus, by (\ref{v=0**}), one obtains
$$\|v_n\|=o(1),$$
a contradiction to (\ref{s(6)}) again.

\ms

\it Case 3. Let $(\a,\b)\in R_3$, i.e., \rm
\bel{}0\le 1-2\a <\b\le{1\over2},\q \a\le{1\over 2}, \qq\m\equiv\m(\a,\b)=2(2\a+\b)-2.\ee
Multiplying the (\ref{s(3)}) by $\l_n^{-1}|\l_n|^{2\b}$, we get
\bel{2.52}i|\l_n|^{-\m+2\b}\|w_n\|^2+\l_n^{-1}|\l_n|^{-\m+2\b}\lan
A^\a v_n,w_n\ran=o(1).\ee
By (\ref{2.40*}) and (\ref{s(4)}),
\bel{2.26}\ba{ll}
\ns\ds|\l_n|^{-1-\m+2\b}\big|\lan A^\a v_n,w_n\ran\big|
\3n\1n&\le|\l_n|^{{\m\over2}+1-4\a}\|A^{\a-{\b\over2}}v_n\|\big(|\l_n|^{-{\m\over2}}\|A^{\b\over2}w_n\|
\big)\\
\ns\ds&\le\big(|\l_n|^{-(2\a-\b)}\|A^{\a-{\b\over2}}v_n\|\big)
\big(|\l_n|^{-{\m\over2}}\|A^{\b\over2}w_n\|\big)=o(1).\ea\ee
Then we obtain from (\ref{2.52}) that
\bel{w=0}|\l_n|^{2-4\a}\|w_n\|^2=o(1),\ee
which implies
\bel{w=0a}\|w_n\|=o(1).\ee
Next, applying bounded operator $A^{\a-{1\over2}}$ to the first
equation in (\ref{2.16a}), we have
\bel{}i\l_n|\l_n|^{-\m}A^\a u_n-|\l_n|^{-\m}A^\a v_n=o(1).\ee
This allows us to rewrite (\ref{2.16c}) as
\bel{2.57}i\l_n|\l_n|^{-\m}w_n+i\l_n|\l_n|^{-\m}A^\a
u_n+|\l_n|^{-\m}A^\b w_n=o(1).\ee
Note that for any $\n\in[0,1]$, ${\n+1\over2}\in[{1\over2},1]$.
Hence, by interpolation, we have
\bel{2.58}|\l_n|^{-\n}\|A^{\n+1\over2}(u_n-A^{\a-1}w_n)\|\le|\l_n|^{-\n}
\|A(u_n-A^{\a-1}w_n)\|^\n\|A^{1\over2}(u_n-A^{\a-1}w_n)\|^{1-\n}=O(1)\ee
due to (\ref{s(7)}) and $\a\le{1\over 2}$.
By taking $\n=1-2\a\in[0,{1\over2}]$, we obtain
\bel{2.58*}\ba{ll}
\ns\ds|\l_n|^{2\a-1}\|A^{1-\a}u_n-w_n\|\3n\1n&=|\l_n|^{-\n}\|A^{\n+1\over2}
(u_n-A^{\a-1}w_n)\|\\
\ns\ds&\le\|\l_n^{-1}(Au_n-A^\a
w_n)\|^{1-2\a}\|A^{1\over2}u_n-A^{\a-{1\over2}}w_n\|^{2\a}=O(1).\ea\ee
Since $\a+\b<1$ in $R_3$ which leads to $\m=2(\b+2\a)-2<2\a$. Hence,
$\m-1<2\a-1$. We now take the inner product of (\ref{2.57}) with
$\l_n^{-1}|\l_n|^{\m}(A^{1-\a}u_n-w_n)$ in $H$,
\bel{2.59}i\lan w_n, A^{1-\a}u_n-w_n\ran+i\|A^{1\over2}u_n\|^2-i\lan
A^\a u_n,w_n\ran+\l_n^{-1}\lan A^\b w_n,
A^{1-\a}u_n-w_n\ran=o(1).\ee
Observe that
$$|\lan w_n,A^{1-\a}u_n-w_n\ran|=\big(|\l_n|^{1-2\a}\|w_n\|\big)\big(|\l_n|^{2\a-1}
\|A^{1-\a}u_n-w_n\|\big)=o(1),$$
due to (\ref{w=0}) and (\ref{2.58*}). It is obvious that the third
term in (\ref{2.59}) is an $o(1)$ because of $\a\le{1\over2}$.
Furthermore, since $1-2\a+\b\in(0,1)$ in $R_3$, we take
$\n=1-2\a+\b$ in (\ref{2.58}) to obtain
$$|\l_n|^{-1+2\a-\b}\|A^{1-\a+{\beta\over2}}(u_n-A^{\a-1}w_n)\|=O(1).
$$
Combining this estimate with (\ref{s(4)}) and the fact $2\b\le1$, we
get
$$\ba{ll}
\ns\ds|\l_n|^{-1}|\lan A^\b w_n,A^{1-\a}u_n-w_n\ran|\\
\ns\ds\le|\l_n|^{-1+2\b}\big(|\l_n|^{-{\m\over2}}\|A^{\b\over2}w_n\|\big)
\big(|\l_n|^{-1+2\a-\b}\|A^{1-2\a+{\b\over2}}(u_n-A^{\a-1}w_n)\|\big)=o(1),\ea$$
i.e., the fourth term in (\ref{2.59}) also converges to zero.
Therefore, we have proved
\bel{}\|A^{1\over2}u_n\|=o(1),\ee
which contradicts (\ref{s(5)}).

\ms

\it Case 4. Let $(\a,\b)\in R_4\cup R_5 \cup S_I$, i.e., \rm
$$0\le\b<\a,\q{1\over2}\le\a,\qq\m={\b\over\a}.$$
By interpolation and (\ref{s(7)}),
\bel{2.31}|\l_n|^{-{\b\over\a}}\|A^\b(v_n+A^{\b-\a}w_n)\|\le
\|\l_n^{-1}A^\a(v_n+A^{\b-\a}w_n)\|^{\b\over\a}\|v_n+A^{\b-\a}w_n\|^{1-{\b\over\a}}
=O(1).\ee
Applying bounded operator $A^{\b-\a}$ to (\ref{2.16c}) leads to
\bel{2.32}i\l_n|\l_n|^{-{\b\over\a}}A^{\b-\a}w_n+|\l_n|^{-{\b\over\a}}A^\b
(v_n+A^{\b-\a}w_n)=o(1).\ee
It follows from (\ref{2.31})--(\ref{2.32}) that
\bel{2.33}|\lambda_n|^{\a-\b\over\a}\|A^{\b-\a}w_n\|=O(1).\ee
Consequently,
\bel{}\ba{ll}
\ns\ds\|w_n\|=\|A^{\a-\b}(A^{\b-\a}w_n)\|\le\|A^{\a-{\b\over2}}(A^{\b-\a}w_n)\|^{\a-\b\over\a-\b/2}
\,\|A^{\b-\a}w_n\|^{\b/2\over\a-\b/2}\\
\ns\ds\qq=\|A^{\b\over2}w_n\|^{2(\a-\b)\over2\a-\b}\,\|A^{\b-\a}w_n\|^{\b\over2\a-\b}\\
\ns\ds\qq=\big(|\l_n|^{-{\b\over2\a}}\|A^{\b\over2}w_n\|\big)^{2(\a-\b)\over2\a-\b}
\big(\l_n^{\a-\b\over\a}\|A^{\b-\a}w_n\|\big)^{\b\over2\a-\b}=o(1).\ea\ee
Here, we have used (\ref{s(4)}) and (\ref{2.33}), and the identity
$$-{\b\over2\a}{2(\a-\b)\over2\a-\b}+
{\a-\b\over\a}{\b\over2\a-\b}=0.$$
Next, note that in region $R_4\cup R_5 \cup S_I$, $1-\a<{1\over 2}$ and $1-2\a+\b<0$
 By applying
$|\l_n|^{\m-1}A^{-{1\over2}}$ to (\ref{2.16a}), we see that
$$\|u_n\|=|\l_n|^{-1}\|v_n\|+o(1)=o(1).$$
Thus, by the boundedness of $\|A^{1\over2}u_n\|$, making use of
interpolation, one gets that $\|A^{1-\a}u_n\|=o(1)$. Moreover, we also have
$\|A^{1-2\a+\b}u_n\|=o(1)$

We take the inner product of (\ref{2.16b}) with
$\l_n^{-1}|\l_n|^{\b\over\a}A^{1-2\a+\b}u_n$ and (\ref{2.16c}) with
$\l_n^{-1}|\l_n|^{\b\over\a}A^{1-\a}u_n$ in H, respectively, to get
the following:
\bel{2.35}i\lan
v_n,A^{1-2\a+\b}u_n\ran+\|\l_n^{-1}A^{1-\a+{\b\over2}}u_n\|^2
-\l_n^{-1}\lan A^\b w_n, A^{1-\a}u_n\ran=o(1),\ee
and
\bel{2.36}i\lan w_n,A^{1-\a}u_n\ran+\l_n^{-1}\lan A^{1\over2}
v_n,A^{1\over2} u_n\ran+\l_n^{-1}\lan A^\b
w_n,A^{1-\a}u_n\ran=o(1).\ee
The first terms in (\ref{2.35}) and (\ref{2.36}) converge to zero,
respectively. We can replace $\l_n^{-1}A^{1\over2}v_n$ in
(\ref{2.36}) by $iA^{1\over2}u_n$ due to (\ref{2.16a}).
Consequently, the sum of (\ref{2.35}) and (\ref{2.36}) yields
$$i\|A^{1\over2}u_n\|^2+\|\l_n^{-1}A^{1-\a+{\b\over2}}u_n\|^2=o(1),$$
which implies
\bel{}\|A^{1\over2}u_n\|=o(1),\ee
a contradiction to (\ref{s(5)}) again. \endpf

\ms

\bf Remark 3.3. \rm In the region $R_2$, $\m=2(2\a-\b)$ stays
unchanged  on the line parallel to the common boundary of $R_2$ and
$R_1$, i.e., the line $\b=2\a-{1\over2}$. It tends to 1 as the
points in $R_2$ get closer to this common boundary.  In the region
$R_3$, the situation is different since the common boundary of $R_3$
and $R_1$ is a single point. In this case, $\m=2(\b+ 2\a)-2$ stays
unchanged on the line parallel to a part of the boundary of $R_3$,
i.e., $\b=-2\a+1$. It tends to 1 as the points in $R_3$ get closer
to the common boundary of $R_3$ and $R_1$. The most interesting case
is the region $R_4\cup R_5$ where $\m={\b\over\a}$ varies on the
line parallel to the common boundary of $R_4$ and $R_1$ but stays
unchanged  on the lines passing the origin. It still tends to 1 as
points in $R_4$ gets closer to the common boundary of $R_4$ and
$R_1$. Moreover, $\m$ is continuous on the region $R_1\cup R_2\cup
R_3\cup R_4$. These observations make us to believe that the orders
of {\it Gevrey} class obtained above are quite reasonable.

\ms

\bf Remark 3.4. \rm The smoothing region given in \cite{MR} does not
include the region $R_5=\{(\a,\b)\bigm|0<\b \le 2\a-1\}$. From the
stability analysis in \cite{HL}, system (\ref{1.2}) is unstable in
this region. However, the instability is caused by the fact that the
origin becomes a spectral point of $\cA_{\a,\b}$, while the Gevrey
class property relies on the behavior of spectrum and resolvent
operator of $\cA_{\a,\b}$ near infinity.

\ms

\section{Asymptotic Behavior of Eigenvalues}
\setcounter{equation}{0}
\setcounter{theorem}{0}

\rm

In this section, we are going to study the asymptotic behavior of
some eigenvalue sequence for the operator $\cA_{\a,\b}$. This will
lead to a proof of Theorem 2.4. Recall that we assume that there
exists a sequence $\m_n$ of eigenvalues of $A$ such that
$$0<\m_1\le\m_2\le\cds,\qq\lim_{n\to\infty}\m_n=\infty.$$
We now present the following lemma.

\ms

\bf Lemma 4.1. \sl Let
\bel{f(l,m)}f(\l,\m)=\l^3+\l^2\m^\b+\l(\m+\m^{2\a})+\m^{\b+1},\qq\forall(\l,\m)\in\dbC\times
\dbR_+. \ee
If the following holds:
\bel{f=0}f(\l_n,\m_n)=0, \ee
then $\l_n$ is an eigenvalue of $\cA_{\a,\b}$.

\ms

\it Proof. \rm For any $\l\in\dbC$, we consider the following
equation for some non-zero $U=(u,v,w)^T\in\cD(\cA_{\a,\b})$ such
that
\bel{}(\l-A_{\a,\b})U=\begin{pmatrix}\l&-I&0\\
A&\l&-A^\a\\0&A^\a&\l+A^\b\end{pmatrix}\begin{pmatrix}u\\ v\\
w\end{pmatrix}=\begin{pmatrix}\l u-v\\ Au+\l v-A^\a w\\
A^\a v+(\l+A^\b)w\end{pmatrix}=0.\ee
Thus,
$$v=\l u,$$
$$w=A^{-\a}(Au+\l v)=A^{1-\a}u+\l A^{-\a}(\l u)=(A^{1-\a}+\l^2A^{-\a})u,$$
and
$$\ba{ll}
\ns\ds0=A^\a(\l u)+(\l+A^\b)(A^{1-\a}+\l^2A^{-\a})u\\
\ns\ds\q=(\l A^\a+\l A^{1-\a}+A^{1+\b-\a}+\l^3A^{-\a}+\l^2 A^{\b-\a})u\\
\ns\ds\q=[\l^3+\l^2A^\b+\l(A+A^{2\a})+A^{\b+1}]A^{-\a}u\equiv
f(\l,A)A^{-\a}u,\ea$$
with $f(\cd\,,\cd)$ given by (\ref{f(l,m)}). Hence, if we take
$u=\f_n$ to be an eigenvector of $A$ corresponding to
$\m_n\in\s(A)$, and let
$$U_n(\l)=\begin{pmatrix}\f_n\\ \l\f_n\\ (\m_n^{1-\a}+\l^2\m_n^{-\a})\f_n\end{pmatrix},$$
then
$$(\l-\cA_{\a,\b})U_n(\l)=\begin{pmatrix}0\\ 0\\ \m_n^{-\a}f(\l,\m_n)\f_n\end{pmatrix}.$$
Therefore, if $\l_n$ is a root of $f(\l,\m_n)=0$, then $\l_n$ is an
eigenvalue of $\cA_{\a,\b}$. \endpf

\ms

Now, for any $n\ge1$, we consider the following equation:
\bel{f=0}f(\l,\m_n)\equiv\l^3+\m_n^\b\l^2+(\m_n^{2\a}+\m_n)\l+\m_n^{\b+1}=0.\ee
Let us denote
\bel{b_n}b_n=\m_n^\b,\q c_n=\m_n^{2\a}+\m_n,\q d_n=\m_n^{\b+1}.\ee
Then (\ref{f=0}) takes the following form:
\bel{3.5}\l^3+b_n\l^2+c_n\l+d_n=0,\ee
with $b_n,c_n,d_n\in\dbR_+$. Let
\bel{pq}p_n=3^2c_n-3b_n^2,\qq q_n=2b_n^3-3^2b_nc_n+3^3d_n.\ee
Define
\bel{D}\ba{ll}
\ns\ds\D_n=\({q_n\over2}\)^2+\({p_n\over3}\)^3=\(b_n^3-{3^2\over2}b_nc_n
+{3^3\over2}d_n\)^2+(3c_n-b_n^2)^3\\
\ns\ds=b_n^6+{3^4\over2^2}b_n^2c_n^2+{3^6\over2^2}d_n^2-3^2b_n^4c_n+3^3b_n^3d_n-{3^5\over2}b_nc_nd_n
+3^3c_n^3-3^3c_n^2b_n^2+3^2c_nb_n^4-b_n^6\\
\ns\ds={3^6\over2^2}d_n^2+3^3b_n^3d_n+3^3c_n^3-{3^5\over2}b_nc_nd_n-{3^3\over2^2}b_n^2c_n^2\\
\ns\ds={3^3\over2^2}\big(3^3d_n^2+2^2b_n^3d_n+2^2c_n^3-2\cd3^2b_nc_nd_n-b_n^2c_n^2\big)\\
\ns\ds={27\over4}\big(27d_n^2+4b_n^3d_n+4c_n^3-18b_nc_nd_n-b_n^2c_n^2\big),\ea\ee
and
\bel{F}\F_{n,\pm}=-{q_n\over2}\pm\sqrt{\D_n}\equiv-{q_n\over2}\pm
\sqrt{\({q_n\over2}\)^2+\({p_n\over3}\)^3}\;.\ee
With the above notations, we have the following result
(\cite{Meserve 1953}).

\ms

\bf Proposition 4.2. (Cardano's Formula). \sl Equation $(\ref{3.5})$
admits three roots which are given by the following:
\bel{root}\l_k={1\over3}\(\F_{n,+}^{1\over3}\o^k+\F_{n,-}^{1\over3}\bar\o^k-b_n\),\q
k=0,1,2, \ee
with $\o=e^{i{2\pi\over3}}\equiv-{1\over2}+i{\sqrt3\over2}$, and for
any $\z=|\z|e^{i\th}$, we define
$\z^{1\over3}=|\z|^{1\over3}e^{i{\th\over3}}$.

\ms

\rm

We note that in the case $\D_n>0$, $\F_{n,\pm}$ are real.
Consequently, the cubic equation (\ref{3.5}) admits a unique real
root, denoted by $\l_{n,0}$ and a pair of complex roots, denoted by
$\l_{n,\pm}$. More precisely, in this case,
\bel{lambdaroot}\left\{\2n\ba{ll}
\ns\ds\l_{n,0}={\F_{n,+}^{1\over3}+\F_{n,-}^{1\over3}-\m_n^\b\over3}\,,\\
\ns\ds\l_{n,\pm}=-{\F_{n,+}^{1\over3}+\F_{n,-}^{1\over3}+\m_n^\b\over6}\pm
i{\sqrt3(\F_{n,+}^{1\over3}-\F_{n,-}^{1\over3})\over6}.\ea\right.\ee
By the definition of $b_n,c_n,d_n$, we have
$$\ba{ll}
\ns\ds\D_n\equiv\D_n(\a,\b)={27\over4}\big(27d_n^2+4b_n^3d_n+4c_n^3-18b_nc_nd_n-b_n^2c_n^2\big)\\
\ns\ds={27\over4}\[27\m_n^{2\b+2}+4\m_n^{4\b+1}+4(\m_n^{2\a}+\m_n)^3
-18\m_n^{2\b+1}(\m_n^{2\a}+\m_n)-\m_n^{2\b}(\m_n^{2\a}+\m_n)^2\]\\
\ns\ds={27\over4}\(27\m_n^{2\b+2}+4\m_n^{4\b+1}+4\m_n^{6\a}+12\m_n^{4\a+1}
+12\m_n^{2\a+2}+4\m_n^3\\
\ns\ds\qq\qq-18\m_n^{2\a+2\b+1}-18\m_n^{2\b+2}-\m_n^{4\a+2\b}-2\m_n^{2\a+2\b+1}-
\m_n^{2\b+2}\)\\
\ns\ds={27\over4}\(8\m_n^{2\b+2}+4\m_n^{4\b+1}+4\m_n^{6\a}+12\m_n^{4\a+1}
+12\m_n^{2\a+2}+4\m_n^3-20\m_n^{2\a+2\b+1}-\m_n^{4\a+2\b}\)\\
\ns\ds=54\m_n^{2\b+2}+27\m_n^{4\b+1}+27\m_n^{6\a}+81\m_n^{4\a+1}
+81\m_n^{2\a+2}+27\m_n^3-135\m_n^{2\a+2\b+1}-{27\over4}\m_n^{4\a+2\b},\ea$$
and
$$\ba{ll}
\ns\ds q_n\equiv q_n(\a,\b)=2b_n^3-3^2b_nc_n+3^3d_n=2\m_n^{3\b}-9\m_n^\b(\m_n^{2\a}+\m_n)+27\m_n^{\b+1}\\
\ns\ds\q\;=2\m_n^{3\b}-9\m_n^{2\a+\b}+18\m_n^{\b+1}.\ea$$
Our first result is about the leading term in $\D_n(\a,\b)$ and in
$q_n(\a,\b)$.

\ms

\bf Lemma 4.3. \sl The following hold:
\bel{D_n}\D_n(\a,\b)=\left\{\2n\ba{ll}
\ns\ds27\m_n^{4\b+1}\big(1+o(1)\big),\qq(\a,\b)\in R_2\cup S_1,\\
\ns\ds216\m_n^3\big(1+o(1)\big),\qq(\a,\b)\in
R_3,\q\a={1\over2},~0\le\b<{1\over2},\\
\ns\ds108\m_n^3\big(1+o(1)\big),\qq(\a,\b)\in R_3\cup S_2,\q0\le\a<{1\over2},~\b={1\over2},\\
\ns\ds27\m_n^3\big(1+o(1)\big),\qq(\a,\b)\in R_3\cup S_2,\q0\le\a,\b<{1\over2},\\
\ns\ds27\m_n^{6\a}\big(1+o(1)\big),\qq(\a,\b)\in R_4\cup R_5\cup
S_I,\ea\right.\ee
and
\bel{q_n}q_n(\a,\b)=\left\{\2n\ba{ll}
\ns\ds2\m_n^{3\b}\big(1+o(1)\big),\qq\qq(\a,\b)\in R_2\cup S_1,\\
\ns\ds9\m_n^{\b+1}\big(1+o(1)\big),\qq(\a,\b)\in
R_3,\q\a={1\over2},~0\le\b<{1\over2},\\
\ns\ds20\m_n^{3\over2}\big(1+o(1)\big),\qq(\a,\b)\in R_3\cup S_2,\q0\le\a<{1\over2},~\b={1\over2},\\
\ns\ds18\m_n^{\b+1}\big(1+o(1)\big),\qq(\a,\b)\in R_3\cup S_2,\q0\le\a,\b<{1\over2},\\
\ns\ds-9\m_n^{2\a+\b}\big(1+o(1)\big),\qq(\a,\b)\in R_4\cup R_5\cup
S_I.\ea\right.\ee

\ms

\rm

\it Proof. \rm For $(\a,\b)\in R_2$, ${1\over2}\le\a<\b\le1$, we
have
$$4\b+1>4\a+2\b\left\{\2n\ba{ll}
\ns\ds\ge2\a+2\b+1\ge2\b+2,\\
\ns\ds>6\a\ge4\a+1\ge2\a+2\ge3.\ea\right.$$
For $(\a,\b)\in R_2\cup S_1$, $0\le\a<{1\over2}<\b$,
$$4\b+1>2\b+2>\left\{\2n\ba{ll}
\ns\ds2\a+2\b+1>4\a+2\b,\\
\ns\ds3>2\a+2>4\a+1>6\a.\ea\right.$$
Hence, in $R_2\cup S_1$,
$$\D_n(\a,\b)=27\m_n^{4\b+1}\big(1+o(1)\big).$$
Also, for $(\a,\b)\in R_2\cup S_1$, $\b>{1\over2}\vee\a$. Thus,
$$3\b>(2\a+\b)\vee(\b+1),$$
which implies that
$$q_n(\a,\b)=2\m_n^{3\b}\big(1+o(1)\big).$$

\ms

For $(\a,\b)\in R_3\cup S_2$, $0\le\a,\b\le{1\over2}$,
$(\a,\b)\ne({1\over2},{1\over2})$. We look at three different cases.

\ms

For $(\a,\b)\in R_3$ with $\a={1\over2}$, and $\b<{1\over2}$, we
have
$$\D_n\({1\over2},\b\)={27\over4}\(32\m_n^3+4\m_n^{4\b+1}
-13\m_n^{2\b+2}\),$$
whose leading term is $216\m_n^3$ since
$$3>2\b+2>4\b+1.$$
Also,
$$q_n\({1\over2},\b\)=2\m_n^{3\b}+9\m_n^{\b+1},$$
whose leading term is $9\m_n^{\b+1}$ since $\b<{1\over2}$.

\ms

For $(\a,\b)\in R_3\cup S_2$ with $\b={1\over2}$,
$0\le\a<{1\over2}$, we have
$$\D_n\(\a,{1\over2}\)={27\over4}\(16\m_n^3+4\m_n^{6\a}+11\m_n^{4\a+1}
-8\m_n^{2\a+2}\),$$
whose leading term is $108\m_n^3$ since
$$3>2\a+2>4\a+1>6\a.$$
Also,
$$q_n\(\a,{1\over2}\)=20\m_n^{3\over2}-9\m_n^{2\a+{1\over2}},$$
whose leading term is $20\m_n^{3\over2}$ since $\a<{1\over2}$.

\ms

Now for $(\a,\b)\in R_3\cup S_2$ with $\a,\b<{1\over2}$, we have
$$3>\left\{\2n\ba{ll}
\ns\ds2\a+2>\left\{\ba{ll}
\ns\ds4\a+1>6\a,\\
\ns\ds2\a+2\b+1>4\a+2\b,\ea\right.\\
\ns\ds4\b+1>2\b+2.\ea\right.$$
Hence, the leading term of $\D_n(\a,\b)$ is $27\m_n^3$. Also,
since
$$\b+1>(3\b)\vee(2\a+\b),$$
the leading term in $q_n(\a,\b)$ is $18\m_n^{\b+1}$.

\ms

Finally, in $R_4\cup R_5\cup S_I$, $0\le\b\vee{1\over2}<\a\le1$, we
have
$$6\a>\left\{\ba{ll}
\ns\ds4\a+1>2\a+2>3,\\
\ns\ds4\a+2\b>\left\{\ba{ll}
\ns\ds2\a+2\b+1>2\b+2,\\
\ns\ds2\a+4\b>4\b+1.\ea\right.\ea\right.$$
Thus, the leading term of $\D_n(\a,\b)$ is $27\m_n^{6\a}$. Also,
since
$$2\a+\b>(3\b)\vee(\b+1),$$
the leading term of $q_n(\a,\b)$ is $-9\m_n^{2\a+\b}$. \endpf

\ms

The following gives the asymptotic behavior of the solutions to
(\ref{f=0}).

\ms

\bf Theorem 4.4. \sl Let the assumption of Theorem 2.4 hold. Let
$n\ge1$ be large enough. Then
\bel{D>0}\D_n(\a,\b)>0,\qq\forall(\a,\b)\in R_2\cup R_3\cup R_4\cup
R_5\cup S_1\cup S_2\cup S_I,\ee
and $(\ref{f=0})$ admits a real root $\l_{n,0}$ and a pair of
conjugate complex roots $\l_{n,\pm}$. Moreover, the following
asymptotic behavior will hold:

\ms

{\rm(i)} For $(\a,\b)\in R_2\cup S_1$,
\bel{eqn3.10}\left\{\2n\ba{l}
\ns\ds\l_{n,0}=-\m_n^\b\big(1+o(1)\big),\\
\ns\ds\l_{n,\pm}=-{1\over2}\m_n^{2\a-\b}\big(1+o(1)\big)\pm
i\m_n^{1\over2}\big(1+o(1)\big).\ea\right.\ee

\ms

{\rm(ii)} For $(\a,\b)\in R_3\cup S_2$,
\bel{eqn3.11}\left\{\2n\ba{l}
\ns\ds\l_{n,0}=\left\{\2n\ba{ll}
\ns\ds-\m_n^\b\big(1+o(1)\big),\qq\qq\a<{1\over2},\\
\ns\ds-{1\over2}\m_n^\b\big(1+o(1)\big),\qq\qq\a={1\over2},\ea\right.\\
\ms\ds\l_{n,\pm}=\left\{\2n\ba{ll}
\ns\ds-{1\over2}\m_n^{2\a+\b-1}\big(1+o(1)\big)\pm
i\m_n^{1\over2}\big(1+o(1)\big),\qq\a,\b<{1\over2},\\
\ns\ds-{1\over4}\m_n^\b\big(1+o(1)\big)\pm
i\sqrt2\m_n^{1\over2}\big(1+o(1)\big),\qq\a={1\over2},\\
\ns\ds-{1\over 4}
\m_n^{2\a - {1\over2}}\big(1+o(1)\big)
\pm\m_n^{1\over2}\big(1+o(1)\big),\qq\b={1\over2}.\ea\right.\ea\right.\ee

\ms

{\rm(iii)} In region $R_4\cup R_5\cup S_I$,
\bel{eqn3.12}\left\{\2n\ba{l}
\ns\ds\l_{n,0}=-\m_n^{1+\b-2\a}\big(1+o(1)\big),\\
\ns\ds\l_{n,\pm}=-{1\over2}\m_n^\b\big(1+o(1)\big)\pm
i\m_n^\a\big(1+o(1)\big).\ea\right.\ee

\ms

\it Proof. \rm By Lemma 4.3, we have (\ref{D>0}). Therefore, the
cubic equation (\ref{f=0}) has one real root and a pair of complex
conjugate roots when $n$ is large enough:
\bel{eqn3.13}\left\{\2n\ba{l}
\ns\ds\l_{n,0}={1\over3}\(\F_{n,+}^{1\over3}+\F_{n,-}^{1\over3}-\m_n^\b\),\\
\ns\ds\l_{n,\pm}=-{1\over6}\(\F_{n,+}^{1\over3}+\F_{n,-}^{1\over3}
+2\m_n^\b\)\pm
i{\sqrt{3}\over6}\(\F_{n,+}^{1\over3}-\F_{n,-}^{1\over3}\).\ea\right.\ee
In what follows, we are going to find the leading terms of the real
and imaginary part of the root expression in (\ref{eqn3.13}).

\ms

\it Case 1: $(\a,\b)\in R_2\cup S_1$, \rm In this case, one has
$$\D_n=27\m_n^{4\b+1}\big(1+o(1)\big),\qq q_n=2\m_n^{3\b}\big(1+\o(1)\big).$$
Thus,
$$\sqrt{\D_n}=3\sqrt3\m_n^{2\b+{1\over2}}\big(1+o(1)\big).$$
Then
$$\ba{ll}
\ns\ds\F_{n,\pm}=-{q_n\over2}\pm\sqrt{\D_n}=-\m_n^{3\b}\big(1+o(1)\big)
\pm3\sqrt3
\m_n^{2\b+{1\over2}}\big(1+o(1)\big)=-\m_n^{3\b}\big(1+o(1)\big),\ea$$
since for $(\a,\b)\in R_2\cup S_1$,
$$3\b>2\b+{1\over2}.$$
Therefore,
\bel{}\ba{ll}
\ns\ds\F_{n,+}^{1\over3}-\F_{n,-}^{1\over3}
={\F_{n,+}-\F_{n,-}\over\F_{n,+}^{2\over3}+\F_{n,+}^{1\over3}
\F_{n,-}^{1\over3}+\F_{n,-}^{2\over3}}
={2\sqrt{\D_n}\over\F_{n,+}^{2\over3}+\F_{n,+}^{1\over3}\F_{n,-}^{1\over3}
+\F_{n,-}^{2\over3}}\\
\ns\ds\qq\qq={6\sqrt3\m_n^{2\b+{1\over2}}\big(1+o(1)\big)\over
3\m_n^{2\b}\big(1+o(1)\big)}=2\sqrt3\m_n^{1\over2}\big(1+o(1)\big),\ea\ee
and
\bel{}\F_{n,+}^{1\over3}+\F_{n,-}^{1\over3}
=-2\m_n^\b\big(1+o(1)\big).\ee
Consequently,
$$\l_{n,0}={1\over3}\(\F_{n,+}^{1\over3}+\F_{n,-}^{1\over3}-\m_n^\b\)
=-\m_n^\b\big(1+o(1)\big).$$
Also,
$$\ba{ll}
\ns\ds\l_{n,\pm}=-{1\over6}\(\F_{n,+}^{1\over3}+\F_{n,-}^{1\over3}
+2\m_n^\b\)\pm
i{\sqrt{3}\over6}\(\F_{n,+}^{1\over3}-\F_{n,-}^{1\over3}\)\\
\ns\ds\qq=-{1\over6}\(\F_{n,+}^{1\over3}+\F_{n,-}^{1\over3}
+2\m_n^\b\)\pm i\m_n^{1\over2}\big(1+o(1)\big).\ea$$
Note that the real part of $\l_{n,\pm}$ above cannot be estimated
using the above argument due to cancelation of the leading term
$\m_n^\b$. Therefore, we take a different approach. To this end, we
denote
$$\L_{n,0}=2\Re\l_{n,\pm}.$$
By the Vieta's formula for the cubic equation (\ref{f=0}), we have
\bel{3.21}-\m_n^\b = \l_{n,0} + \l_{n,+} + \l_{n,-} = \l_{n,0} +
2\mbox{Re}\l_{n,\pm}=\l_{n,0}+\L_{n,0}.\ee
Therefore, $\l_{n,0}=-\L_{n,0}-\m_n^\b$ satisfies (\ref{f=0}), i.e.,
$$\ba{ll}
\ns\ds0=\l_{n,0}^3+\m_n^\b\l_{n,0}^2+(\m_n^{2\a}+\m_n)\l_{n,0}+\m_n^{\b+1}\\
\ns\ds\q=-(\L_{n,0}+\m_n^\b)^3+\m_n^\b(\L_{n,0}+\m_n^\b)^2-(\m_n^{2\a}+\m_n)(\L_{n,0}+\m_n^\b)+\m_n^{\b+1}\\
\ns\ds\q=-\L_{n,0}^3-3\m_n^\b\L_{n,0}^2-3\m_n^{2\b}\L_{n,0}-\m_n^{3\b}+\m_n^\b\L_{n,0}^2+2\m_n^{2\b}
\L_{n,0}+\m_n^{3\b}\\
\ns\ds\qq-(\m_n^{2\a}+\m_n)\L_{n,0}-\m_n^{2\a+\b}-\m_n^{\b+1}+\m_n^{\b+1}\\
\ns\ds\q=-\L_{n,0}^3-2\m_n^\b\L_{n,0}^2-(\m_n^{2\a}+\m_n^{2\b}+\m_n)\L_{n,0}-\m_n^{2\a+\b}.\ea$$
This means that $\L_{n,0}$ is a real solution to the following new
cubic equation
\bel{3.23}\L^3+2\m_n^\b\L^2+\big(\m_n^{2\a}+\m_n^{2\b}+\m_n
\big)\L+\m_n^{2\a+\b}=0. \ee
Next, by defining
$$\L_{n,\pm}=-\m_n^\b-\l_{n,\pm},$$
and by the fact that $\l_{n,\pm}$ are the roots of (\ref{f=0}),
using the same argument as above, we see that $\L_{n,\pm}$ is a pair
of conjugate complex roots of (\ref{3.23}). Now, we rewrite
(\ref{3.21}) as follows:
$$\L_{n,0} = -\l_{n,0} - \m_n^\b = 2\mbox{Re}\l_{n,\pm}.$$
Since the leading term of $\mbox{Re}\l_{n,\pm}$ is $o(\m_n^\b)$, the
complex roots of equation (\ref{3.23}) satisfy
$$\L_{n,\pm}=-\m_n^\b-\l_{n,\pm}=-\m_n^\b(1+o(1))\mp i\m_n^{1\over2}(1+o(1)).$$
Then by the Vieta's formula for equation (\ref{3.23}), one has
(noting $\b>{1\over2}$)
$$-\m_n^{2\a+\b}=\L_{n,0}\L_{n,+}\L_{n,-}=\L_{n,0}\m_n^{2\b}(1+o(1)).$$
Therefore,
$$\Re\l_{n,\pm}={1\over2}\L_{n,0}=-{1\over2}\m_n ^{2\a-\b}(1+o(1)).$$

\ms

\it Case 2: \rm $(\a,\b)\in R_4\cup R_5\cup S_I$. In this case,
$$\D_n=27\m_n^{6\a}\big(1+o(1)\big),\qq q_n=-9\m_n^{2\a+\b}\big(1+o(1)\big).$$
Then,
$$\sqrt{\D_n }=3\sqrt3\m_n^{3\a}\big(1+o(1)\big).$$
This leads to
$$\ba{ll}
\ns\ds\F_{n,\pm}=-{q_n\over2}\pm\sqrt{\D_n}
={9\over2}\m_n^{2\a+\b}\big(1+o(1)\big)
\pm3\sqrt3\mu_n^{3\a}\big(1+o(1)\big)=\pm3\sqrt3\m_n^{3\a}\big(1+o(1)\big),\ea$$
since for the current case, $\b<\a$ which implies
$$3\a>2\a+\b.$$
Hence,
\bel{}\F_{n,+}^{1\over3}-\F_{n,-}^{1\over3}=2\F_{n,+}^{1\over3}=
2\sqrt3\m_n^\a\big(1+o(1)\big),\ee
and
\bel{}\ba{ll}
\ns\ds\F_{n,+}^{1\over3}+\F_{n,-}^{1\over3}
={\F_{n,+}+\F_{n,-}\over\F_{n,+}^{2\over3}-\F_{n,+}^{1\over3}
\F_{n,-}^{1\over3}+\F_{n,-}^{2\over3}}
={-q_n\over\F_{n,+}^{2\over3}-\F_{n,+}^{1\over3}\F_{n,-}^{1\over3}
+\F_{n,-}^{2\over3}}\\
\ns\ds\qq\qq\qq={9\m_n^{2\a+\b}\big(1+o(1)\big)\over
9\m_n^{2\a}\big(1+o(1)\big)}=\m_n^\b\big(1+o(1)\big).\ea\ee
Consequently,
$$\ba{ll}
\ns\ds\l_{n,\pm}=-{1\over6}\(\F_{n,+}^{1\over3}+\F_{n,-}^{1\over3}
+2\m_n^\b\)\pm
i{\sqrt{3}\over6}\(\F_{n,+}^{1\over3}-\F_{n,-}^{1\over3}\)
=-{1\over2}\m_n^\b\big(1+o(1)\big)\pm i\m_n^\a\big(1+o(1)\big).\ea$$
For the real root $\l_{n,0}$, we will have cancelation of the
leading term $\m_n^\b$. Therefore, we may let
$$\l_{n,0}=c\m_n^\xi\big(1+o(1)\big),$$
for some $c\in\dbR$ and $0<\xi<\b$. Then by Vieta's formula, noting
$\b<\a$,
$$-\m_n^{\b+1}=\l_{n,0}\l_{n,+}\l_{n,-}=c\m_n^\xi\({1\over4}\m_n^{2\b}
+\m_n^{2\a}\)\big(1+o(1)\big)=c\m^{2\a+\xi}\big(1+o(1)\big).$$
Consequently, it is necessary that
$$c=-1,\qq\xi=\b+1-2\a.$$

\it Case 3: \rm $(\a,\b)\in R_3\cup S_2$. We will consider three
subcases.

\ms

\it Subcase 1: \rm $(\a,\b)\in R_3\cup S_2$, $0\le\a,\b<{1\over2}$.
\rm In this case,
$$\D_n=27\m_n^3\big(1+o(1)\big),\qq q_n=18\m_n^{\b+1}\big(1+o(1)\b).$$
Then,
$$\sqrt{\D_n}=3\sqrt3\m_n^{3\over2}\big(1+o(1)\big).$$
This further gives
$$\ba{ll}
\ns\ds\F_{n,\pm}=-{q_n\over2}\pm\sqrt{\D_n}
=-9\m_n^{\b+1}\big(1+o(1)\big)\pm3\sqrt3
\m_n^{3\over2}\big(1+o(1)\big)=\pm3\sqrt3
\m_n^{3\over2}\big(1+o(1)\big),\ea$$
since for the current case,
$${3\over2}>\b+1.$$
Then
\bel{}\F_{n,+}^{1\over3}-\F_{n,-}^{1\over3}
=2\sqrt3\m_n^{1\over2}\big(1+o(1)\big),\ee
and
\bel{}\ba{ll}
\ns\ds\F_{n,+}^{1\over3}+\F_{n,-}^{1\over3}
={\F_{n,+}+\F_{n,-}\over\F_{n,+}^{2\over3}-\F_{n,+}^{1\over3}
\F_{n,-}^{1\over3}+\F_{n,-}^{2\over3}}
={-q_n\over\F_{n,+}^{2\over3}-\F_{n,+}^{1\over3}\F_{n,-}^{1\over3}
+\F_{n,-}^{2\over3}}\\
\ns\ds\qq\qq\qq={-18\m_n^{\b+1}\big(1+o(1)\big)\over
9\m_n\big(1+o(1)\big)}=-2\m_n^\b\big(1+o(1)\big).\ea\ee
Consequently,
$$\l_{n,0}={1\over3}\(\F_{n,+}^{1\over3}+\F_{n,-}^{1\over3}-\m_n^\b\)
=-\m_n^\b\big(1+o(1)\big),$$
and
$$\ba{ll}
\ns\ds\l_{n,\pm}=-{1\over6}\(\F_{n,+}^{1\over3}+\F_{n,-}^{1\over3}
+2\m_n^\b\)\pm
i{\sqrt{3}\over6}\(\F_{n,+}^{1\over3}-\F_{n,-}^{1\over3}\)\\
\ns\ds\qq =-{1\over6}\(\F_{n,+}^{1\over3}+\F_{n,-}^{1\over3}
+2\m_n^\b\)\pm i\m_n^{1\over2}\big(1+o(1)\big).\ea$$
From the above, we see that the leading terms in $\Re\l_{n,\pm}$ are
canceled. Thus,
$$\L_{n,0}\equiv2\Re\l_{n,\pm}=o(\m_n^\b).$$
The same as Case 1, $\L_{n,0}$ is a real root of cubic equation
(\ref{3.23}), and
$$\L_{n,\pm}=-\m_n^\b-\l_{n,\pm}=-\m_n^\b-\Re\l_{n,\pm}
\pm i\m_n^{1\over2}\big(1+o(1)\big)=-\m_n^\b\big(1+o(1)\big)\pm
i\m_n^{1\over2}\big(1+o(1)\big)$$
are the pair of conjugate complex roots of (\ref{3.23}). Further, by
Vieta's formula for the equation (\ref{3.23}), one obtains
$$\ba{ll}
\ns\ds-\m_n^{2\a+\b}=\L_{n,0}\L_{n,+}\L_{n,-}=\L_{n,0}\[\(\Re\L_{n,\pm}\)^2+\(\Im\L_{n,\pm}\)^2\]\\
\ns\ds\qq\qq=\L_{n,0}\[\m_n^{2\b}\big(1+o(1)\big)+\m_n\big(1+o(1)\big)\]
=\L_{n,0}\m_n\big(1+o(1)\big),\ea$$
since $\b<{1\over 2}$. Consequently,
$$\L_{n,0}=-\m_n^{2\a+\b-1}\big(1+o(1)\big).$$
Hence,
$$\Re\l_{n,\pm}={1\over2}\L_{n,0}=-{1\over2}\m_n^{2\a+\b-1}\big(1+o(1)\big),$$
proving our claim.

\ms

\it Subcase 2: \rm $(\a,\b)\in R_3$, $\a={1\over2}$, $\b<{1\over2}$.
For this case,
$$\D_n=216\m_n^3\big(1+o(1)\big),\qq q_n=9\m_n^{\b+1}\big(1+o(1)\big).$$
Then
$$\ba{ll}
\ns\ds\F_{n,\pm}=-{q_n\over2}\pm\sqrt{\D_n}=-{9\over2}\m_n^{\b+1}
\big(1+o(1)\big)
\pm6\sqrt6\m_n^{3\over2}\big(1+o(1)\big)=\pm6\sqrt6\m_n^{3\over2}\big(1+o(1)\big),\ea$$
since
$${3\over2}>\b+1.$$
We have
\bel{}\F_{n,+}^{1\over3}-\F_{n,-}^{1\over3}
=2\sqrt6\m_n^{1\over2}\big(1+o(1)\big),\ee
and
\bel{}\ba{ll}
\ns\ds\F_{n,+}^{1\over3}+\F_{n,-}^{1\over3}
={\F_{n,+}+\F_{n,-}\over\F_{n,+}^{2\over3}-\F_{n,+}^{1\over3}
\F_{n,-}^{1\over3}+\F_{n,-}^{2\over3}}
={-q_n\over\F_{n,+}^{2\over3}-\F_{n,+}^{1\over3}\F_{n,-}^{1\over3}
+\F_{n,-}^{2\over3}}\\
\ns\ds\qq\qq\qq={-9\m_n^{\b+1}\big(1+o(1)\big)\over
18\m_n\big(1+o(1)\big)}=-{1\over2}\m_n^\b\big(1+o(1)\big).\ea\ee
Consequently,
$$\l_{n,0}={1\over3}\(\F_{n,+}^{1\over3}+\F_{n,-}^{1\over3}-\m_n^\b\)
=-{1\over2}\m_n^\b\big(1+o(1)\big),$$
and
$$\ba{ll}
\ns\ds\l_{n,\pm}=-{1\over6}\(\F_{n,+}^{1\over3}+\F_{n,-}^{1\over3}
+2\m_n^\b\)\pm
i{\sqrt{3}\over6}\(\F_{n,+}^{1\over3}-\F_{n,-}^{1\over3}\)\\
\ns\ds\qq =-{1\over4}\m_n^\b\pm
i\sqrt2\m_n^{1\over2}\big(1+o(1)\big).\ea$$

\it Subcase 3: \rm $(\a,\b)\in R_3\cup S_2$, $\b={1\over2}$ and
$\a<{1\over2}$. In this case,
$$\D_n=108\m_n^3\big(1+o(1)\big),\qq
q_n=20\m_n^{3\over2}\big(1+o(1)\big).$$
Then
$$\F_{n,\pm}=-{q_n\over2}\pm\sqrt{\D_n}=\big(-10\pm6\sqrt3\big)\m_n^{3\over2}\big(1+o(1)\big)
=(-1\pm\sqrt3)^3\m_n^{3\over2}\big(1+o(1)\big).$$
Hence,
\bel{}\F_{n,+}^{1\over3}-\F_{n,-}^{1\over3}
=\[(-1+\sqrt3)-(-1-\sqrt3)\]\m_n^{1\over2}
\big(1+o(1)\big)=2\sqrt3\m_n^{1\over2}\big(1+o(1)\big),\ee
and
\bel{}\F_{n,+}^{1\over3}+\F_{n,-}^{1\over3}
=\[(-1+\sqrt3)+(-1-\sqrt3)\]\m_n^{1\over2} \big(1+o(1)\big)=-2
\m_n^{1\over2}\big(1+o(1)\big).\ee
Consequently,
$$\l_{n,0}={1\over3}\(\F_{n,+}^{1\over3}+\F_{n,-}^{1\over3}-\m_n^{1\over2})
=-\m_n^{1\over2}\big(1+o(1)\big),$$
and
$$\ba{ll}
\ns\ds\l_{n,\pm}=-{1\over6}\(\F_{n,+}^{1\over3}+\F_{n,-}^{1\over3}
+2\m_n^{1\over2}\)\pm
i{\sqrt{3}\over6}\(\F_{n,+}^{1\over3}-\F_{n,-}^{1\over3}\)\\
\ns\ds\qq=-{1\over6}\(\F_{n,+}^{1\over3}+\F_{n,-}^{1\over3}
+2\m_n^{1\over2}\)\pm i\m_n^{1\over2}\big(1+o(1)\big).\ea$$
Once again, we see that the leading terms in $\Re\l_{n,\pm}$ are
canceled out. Therefore,
$$ \Re\l_{n,\pm}=o(\m_n^{1\over2}).$$
Mimicking Case 1, we know that
$$\left\{\2n\ba{ll}
\ns\ds\L_{n,0}=2\Re\l_{n,\pm},\\
\ns\ds\L_{n,\pm}=-\m_n^{1\over 2}-\l_{n,\pm}=-\m_n^{1\over
2}(1+o(1)) \mp\m_n^{1\over 2}\big(1+o(1)\big)\ea\right.$$
are three roots of the cubic equation (\ref{3.23}). Then by the
Vieta's formula,
$$ -\m_n^{2\a + {1\over 2}}= \L_{n,0}\L_{n,+}\L_{n,-}=\L_{n,0}[2\m_n+(1+o(1))].$$
Therefore,
$$ \L_{n,0}= -{1\over 2}\m_n^{2\a - {1\over 2}}(1+o(1)), $$
i.e.,
$$ \mbox{Re}\l_{n,\pm} = {1\over 2}\L_{n,0} = -{1\over 4}\m_n^{2\a - {1\over 2}}(1+o(1)).$$
This completes the proof. \endpf

\ms

We see easily that (\ref{2.19})--({2.20}) follows from Theorem 4.4.
Therefore, proof of Theorem 2.4 follows immediately.

\ms

\bf Remark 4.5. \rm In our previous paper \cite{HL}, a complete
stability analysis for system (\ref{1.1}) was presented. The
asymptotic expressions of eigenvalues in $\l_{n,0}$ and $\l_{n,\pm}$
for $(\a,\b)\in S_1\cup S_2\cup S_I$ were derived by plugging the
Taylor series expansion of $\Phi_{n,\pm}^{1\over3}$ into
(\ref{lambdaroot}). Due to the cancelation of leading term and other
terms, this method became cumbersome in finding an explicit ordering
of the power terms of $\m_n$ in each region. A number of subregions
were further introduced, but the argument there were not clear and
satisfactory. The idea used in the current paper is much better and
it enable us to present a complete analysis of the asymptotic
behavior of the eigenvalues.

\ms

We now present an interesting corollary of Theorem 4.4, which also
gives us an impression that the index $\m(\a,\b)$ is sharp.

\ms

\bf Corollary 4.6. \sl Under the assumption of Theorem 2.4, the
following holds:
\bel{4.31}\limsup_{\l\in\dbR,|\l|\to\infty}|\l|^{\m(\a,\b)}
\|(i\l-\cA_{\a,\b})^{-1}\|\ge2,\qq\forall(\a,\b)\in[0,1]\times[0,1]\setminus
R_1.\ee

\ms

\it Proof. \rm First of all, we claim that if $\l\in\dbC$, with
$\Re\l\ne0$, is an eigenvalue of $\cA$ which is a densely defined
closed operator on some Hilbert space $\cH$ such that
$(i\Im\l-\cA)^{-1}$ exists. Then
\bel{4.1}|\Re\l|\,\|(i\Im\l-\cA)^{-1}\|\ge1.\ee
In fact, there exists an $x\in\cD(\cA)$ with $\|x\|=1$ such that
$$\cA x=(\m+i\n)x.$$
Hence,
$$\m(i\n-\cA)^{-1}x=-x.$$
Thus (\ref{4.1}) follows. Now, from Theorem 4.4, we know that
$\cA_{\a,\b}$ has a sequence of conjugate complex eigenvalues of the
following form:
$$\l_{n,\pm}=-a\m_n^\xi\big(1+o(1)\big)\pm
ib\m_n^\eta\big(1+o(1)\big),$$
for some real constants $a,b,\eta>0$ and $\xi\ge0$. Let
$\l=b\m_n^\eta$. By (\ref{4.1}), we have
$$\ba{ll}
\ns\ds1\le
a\m_n^\xi\big(1+o(1)\big)\big\|\big(ib\m_n^\eta-\cA_{\a,\b}\big)^{-1}\big\|
={a\over
b^{\xi\over\eta}}|\l|^{\xi\over\eta}\big\|\big(i\l-\cA)^{-1}\big\|\big(1+o(1)\big).\ea$$
Hence,
$$|\l|^{\xi\over\eta}\big\|(i\l-\cA_{\a,\b})^{-1}\big\|\ge{b^{\xi/\eta}\over
a}\big(1+o(1)\big).$$
Now, we look at different regions.

\ms

In region $R_2\cup S_1$,
$$\l_{n,\pm}=-{1\over2}\m_n^{2\a-\b}\big(1+o(1)\big)\pm
i\m_n^{1\over2}\big(1+o(1)\big).$$
Thus,
$${\xi\over\eta}=2(2\a-\b),\qq{b^{\xi/\eta}\over a}=2,$$
which leads to
$$|\l|^{2(2\a-\b)}\big\|(i\l-\cA_{\a,\b})^{-1}\big\|\ge2\big(1+o(1)\big).$$

\ms

In region $R_3\cup S_2$,
\bel{}\l_{n,\pm}=\left\{\2n\ba{ll}
\ns\ds-{1\over2}\m_n^{2\a+\b-1}\big(1+o(1)\big)\pm
i\m_n^{1\over2}\big(1+o(1)\big),\qq\a,\b<{1\over2},\\
\ns\ds-{1\over4}\m_n^\b\big(1+o(1)\big)\pm
i\sqrt2\m_n^{1\over2}\big(1+o(1)\big),\qq\a={1\over2},\\
\ns\ds-{1\over 4}\m_n^{2\a - {1\over2}}\big(1+o(1)\big)
\pm\m_n^{1\over2}\big(1+o(1)\big),\qq\b={1\over2}.\ea\right.\ee
Hence,
$$\left\{\2n\ba{ll}
\ns\ds{\xi\over\eta}=2(2\a+\b-1),\qq{b^{\xi/\eta}\over
a}=2,\qq\a,\b<{1\over2},\\
\ns\ds{\xi\over\eta}=2\b,\qq{b^{\xi/\eta}\over
a}=2^{\b+2},\qq\a={1\over2},\\
\ns\ds{\xi\over\eta}=4\a-1,\qq{b^{\xi/\eta}\over
a}=4,\qq\b={1\over2}.\ea\right.$$
Consequently,
$$|\l|^{2(2\a+\b-1)}\big\|(i\l-\cA_{\a,\b})^{-1}\big\|\ge2.$$

\ms

Finally, in $R_4\cup R_5\cup S_I$,
$$\l_{n,\pm}=-{1\over2}\m_n^\b\big(1+o(1)\big)\pm
i\m_n^\a\big(1+o(1)\big).$$
Thus,
$${\xi\over\eta}={\b\over\a},\qq{b^{\xi/\eta}\over
a}=2,$$
and we again have
$$|\l|^{\b\over\a}\big\|(i\l-\cA_{\a\,\b})^{-1}\big\|\ge2\big(1+o(1)\big).$$
Combining the above, we see that for any given
$(\a,\b)\in[0,1]\times[0,1]$,
$$|\l|^{\m(\a,\b)}\big\|(i\l-\cA_{\a,\b})^{-1}\big\|\ge2\big(1+o(1)\big),$$
with $\m(\a,\b)$ given by (\ref{mu}). Hence, (\ref{4.31}) follows.
\endpf

\ms

To conclude this paper, we point out that with the complete
stability and regularity results for system (\ref{1.1}), we should
be able to consider the more general system (\ref{1.7}) when the
operators $B$ and $D$ are equivalent (in a certain sense) to $A^\a$
and $A^\b$, respectively. Such a general setting will allow
differential operators to have different boundary conditions.
Relevant results will be addressed in a forthcoming paper.

\bs

\no {\bf Acknowledgement}. The authors would like to thank Zhaobin
Kuang, a graduate student at the University of Minnesota Duluth, for
his idea of using the Vieta's formula in the proof of Theorem 4.4.

\end{document}